\documentclass[12pt]{article}

\usepackage{amsmath,amsthm,amsfonts,graphicx,tikz-cd,pgfplots,listings}

\usepackage{geometry,titlesec,hyperref,xhfill,setspace,float,fancyhdr,cite}

\usetikzlibrary{positioning,calc}
\usepgfplotslibrary{polar}
\usepgflibrary{shapes.geometric}
\usepgfplotslibrary{fillbetween}

\pgfplotsset{my style/.append style={axis x line=middle, axis y line=
middle, xlabel={$x$}, ylabel={$y$}, axis equal }}

\theoremstyle{definition}
\newtheorem{definition}{Definition}[section]
\newtheorem{example}{Example}[section]

\theoremstyle{plain}
\newtheorem{theorem}{Theorem}[section]
\newtheorem{lemma}{Lemma}[section]
\newtheorem{corollary}{Corollary}[section]

\newcommand{\underscore}{\underline{\hspace{2mm}}}
\newcommand{\Z}{\mathbb{Z}}

\newcommand{\C}{\mathbb{C}}
\newcommand{\R}{\mathbb{R}}

\newcommand{\Aut}{\text{Aut}}
\newcommand{\M}{\text{M}}
\newcommand{\SL}{\text{SL}}
\newcommand{\GL}{\text{GL}}
\newcommand{\SO}{\text{SO}}
\newcommand{\PGL}{\text{PGL}}
\newcommand{\PSL}{\text{PSL}}
\newcommand{\PSO}{\text{PSO}}

\newcommand{\Hom}{\text{Hom}}
\newcommand{\Out}{\text{Out}}
\newcommand{\RP}{\mathbb{R}P}
\newcommand{\slf}{\mathfrak{sl}}

\newcommand{\pso}{\mathfrak{pso}}

\newcommand{\g}{\mathfrak{g}}
\newcommand{\vv}{\mathfrak{v}}

\newcommand{\Ad}{\text{Ad }}

\newcommand{\lra}{\longrightarrow}

\newcommand{\tr}{\text{tr }}

\newcommand{\ZG}{\mathbb{Z}\Gamma}
\newcommand{\del}{\partial}
\newcommand{\res}{\text{res}}

\title{Projective Rigidity of Once-Punctured Torus Bundles via the Twisted Alexander Polynomial}
\author{Charles Daly}

\date{
	\today
}

\begin{document}
	\maketitle
	
	\begin{abstract}
		In this paper we provide a means of certifying infinitesimal projective rigidity relative to the cusp for hyperbolic once-punctured torus bundles in terms of twisted Alexander polynomials of representations associated to the holonomy.  We also relate this polynomial to an induced action on the tangent space of the character variety of the free group of rank 2 into $\PGL(4,\R)$ that arises from the holonomy of a hyperbolic once-punctured torus bundle.  We prove the induced action on the tangent space of the character variety is the same as the group theoretic action that arises in the Lyndon–Hochschild–Serre spectral sequence on cohomology.
	\end{abstract}

\begin{section}{Overview of the Paper}\label{secoverview}
This paper is an extension of a work by the author to the context of hyperbolic once-punctured torus bundles.  In Daly \cite{Daly}, we provided several cohomology calculations related to the figure-eight knot complement and its holonomy into $\PSO(3,1) \subset \PGL(4,\R)$.  These calculations were motivated to study both the projective deformations of the hyperbolic structure of the figure-eight knot complement and probe its character variety into $\PGL(4,\R)$.  This paper carries out many of those calculations in the greater context of once-punctured torus bundles and relates certain algebraic invariants of this representation to the so-called \emph{twisted Alexander polynomial}.  We provide computational means of verifying infinitesimal projective rigidity rel cusp in terms of roots of this polynomial if we are given the holonomy representation into $\PSO(3,1)$.  These findings are summarized in Theorem \ref{thm2} and Theorem \ref{thm3}.  \\
\\
Section \ref{section-preliminaries} begins with an overview some algebra, specifically twisted group cohomology and the twisted Alexander polynomial of a representation.  We provide examples of some calculations of the twisted Alexander polynomial.  In addition we motivate this study by projective geometry and deformations of the hyperbolic structure of once-punctured torus bundles.  We recall notions of infinitesimal projective rigidity and Weil's Infinitesimal Rigidity theorem.  Section \ref{sectioncc} addresses hyperbolic once-punctured torus bundles.  We provide a bit of background about their group structure, holonomy representations via the trace relations, and provide some calculations of their cohomology twisted by the holonomy representation.  In this section we reinterpret a group action associated to the splitting of the once-punctured torus bundle's fundamental group in terms of deformations of representations of its fundamental group.  This reinterpretation in Theorem \ref{thm1p5} gives geometric meaning to the $G/N$-action on cohomology afforded by Lyndon–Hochschild–Serre spectral sequence where $N$ is a normal subgroup of a group $G$.  In Section \ref{secCalcandExs} we provide some calculations and examples to verify infinitesimal projective rigidity rel cusp for two hyperbolic once-punctured torus bundles.  These calculations can be carried out in far greater generality with the aid of interval arithmetic, however we provide a basic framework and carry out the calculations in particularly nice examples.  These calculations are carried out through Mathematica \cite{Mathematica} and companion programs can be found on the author's webpage \href{https://www.math.brown.edu/cdaly2/research.html}{Daly}.
\end{section}

\section*{Acknowledgements}
I would like to thank several folks for their help writing this paper.  Richard Schwartz has been a constant source of motivation and inspiration throughout this work.  He was the one to originally suggest trying to carry out these cohomology calculations in greater generality from my work on the figure-eight knot complement.  Bill Goldman provided very helpful conversations regarding the outer automorphism action on character varieties and provided me with several helpful references.  I would also like to thank both Michael Heusener and Sam Ballas for their insights into these calculations.  They have provided very important insights and avenues to explore in this work.  I extend my deep gratitude to Joan Porti who hosted me at the Universitat Aut\`{o}noma de Barcelona where this work began to develop a tangible direction.  There he patiently described what the twisted Alexander polynomial is and related my work at the time to it.

\begin{section}{Preliminaries}\label{section-preliminaries}
In this section we provide the necessary preliminaries that we will use to provide a sufficient condition for infinitesimal projective rigidity relative to the cusp of a hyperbolic once-punctured torus bundle.  The first section addresses some basics of twisted group cohomology.  More detailed accounts of these topics can be found in works by Brown and L\"{o}h  \cite{Brown1982Cohomology}, \cite{claraloh}.  The next section provides a brief overview of the twisted Alexander polynomial.  This account is largely taken from the works of Kitano and Wada \cite{Kitano}, \cite{WADA1994241}.  The final part of the this section provides some definitions and motivation for considering the problem at hand regarding infinitesimal projective rigidity rel cusp.

\begin{subsection}{Twisted Cohomology}\label{subsec-twistedcohom}
We begin with a brief overview of twisted group cohomology.  Let $R$ be a ring, and $M$ be an $R$-module.  Recall that an $R$-module $F$ is called \emph{free} if it admits a basis, i.e. a subset $\beta \subset M$ for which every non-zero element $m \in M$ admits a \emph{unique} decomposition as a finite $R$-linear combination of elements of $\beta$.  An $R$-module is called $\emph{projective}$ if every short exact sequence of $R$-modules
\begin{equation*}
0 \lra K \lra E \lra P \lra 0
\end{equation*}
is split-exact, i.e. the map $E \lra P$ admits an $R$-module section.  As a consequence of the definitions, every free module is projective as the existence of a basis gives rise to sections simply by pulling back the basis elements arbitrarily and extending the map linearly.  \\
\\
A (left) \emph{resolution} of $M$ is a long exact sequence of $R$-modules
\begin{equation*}
\hdots \lra N_{2} \lra N_{1} \lra N_{0} \lra M \lra 0
\end{equation*}
We say that a resolution of $M$ is \emph{free}, (respectively projectively), if each such $N_{i}$ is a free, (respectively projective), $R$-module.  Note that the existence of a free-resolution comes immediately by picking a presentation of $M$ and carrying out the extension inductively step-by-step.  \\
\\
Let $\Gamma$ be a group.  Typically speaking we are interested in discrete, finitely presented groups, however the following constructions hold in far greater generality.  The \emph{group ring} $\ZG$ is a (typically non-commutative) ring given by all \emph{finitely supported} $\Z$-linear combinations of $\Gamma$ where multiplication is carried out via
\begin{equation*}
\left(\sum_{g \in G} a_{g} g\right)\left(\sum_{h \in G} b_{h} h\right) = \sum_{g,h \in G} (a_{g}b_{h})(gh)
\end{equation*}
For a given group $\Gamma$, a representation of $\Gamma$ into a finite-dimensional vector space $V$ gives rise to a $\ZG$-module, and conversely, a $\ZG$-module $V$, induces a $\Gamma$-action.  For the context of this paper, we will typically be interested in $\ZG$-modules and $\Z$-resolutions of $\ZG$-modules where $\Z$ is given the trivial $\Gamma$-action.  \\
\\
We now have the necessary terminology to define the twisted group cohomology.  Let $V$ be a finite-dimensional representation of $\Gamma$, i.e. a $\ZG$-module.  Let $P$ denote a projective $\ZG$-resolution of $\Z$.  
\begin{equation*}
\hdots \lra P_{2} \xrightarrow{\del_{2}} P_{1} \xrightarrow{\del_{1}} P_{0} \xrightarrow{\del_{0}} \Z \rightarrow 0
\end{equation*}
We apply the $\Hom_{\ZG}\left( \bullet, V\right)$-functor to the above long-exact sequence to obtain a chain-complex $C^{n}(\Gamma; V) := \Hom_{\ZG}\left(P_{n},V\right)$.  
\begin{equation*}
0 \rightarrow \Hom_{\ZG}(\Z,V) \xrightarrow{\delta^{0}}  C^{0}(\Gamma; V) \xrightarrow{\delta^{1}} C^{1}(\Gamma; V)  \xrightarrow{\delta^{2}} C^{2}(\Gamma; V) \xrightarrow{\delta^{3}} \hdots
\end{equation*}
where $\delta^{n+1}: C^{n}(\Gamma; V) \lra C^{n+1}(\Gamma; V)$ acts via $(\delta^{n+1}f_{n})(p_{n+1}):=f_{n}\left(\del_{n+1}p_{n+1}\right)$ for all $f_{n} \in C^{n}(\Gamma; V)$ and $p_{n+1} \in P_{n+1}$.  We define the \emph{twisted group cohomology of the representation} $V$ via the cochain complex above.  Specifically, $H^{n}(\Gamma; V):= H^{n}\left(C^{*}(G;V)\right)$.  \\
\\
We remark our main interest will be in the \emph{first twisted cohomology group}, $H^{1}(G; V)$.  This is in large part because under ideal circumstances $H^{1}(\Gamma; \g_{\Ad \rho})$ parametrizes the \emph{infinitesimal deformations} of a fixed representation $\rho: \Gamma \lra G$, up to conjugation, where $\g_{\Ad \rho}$ is $\Gamma$-module associated to the original representation $\rho$ post-composed with the adjoint representation, $\Ad: G \lra \Aut(\g)$.  \\
\\
Expanding upon the remarks from the previous paragraph, let $\rho: \Gamma \lra \GL(V)$ be a representation of $\Gamma$ into a finite dimensional vector space $V$.  A \emph{crossed homomorphism} is a map $f: \Gamma \lra V$ satisfying $f(gh) = f(g) + gf(h)$ for all $g,h \in \Gamma$ where here $g$ is acting on $f(h)$ via the representation $\rho$.  A \emph{coboundary} is a crossed homomorphism $f$ of the form $f(g) = (1-g)v$ for some fixed $v \in V$ and for all $g \in G$.  Loosely speaking, Weil's Infinitesimal Rigidity Theorem states that crossed homomorphisms parametrize the tangent space of the \emph{representation variety} $\Hom(\Gamma,G)$ at the point $\rho$, whereas coboundaries parametrize the tangent vectors at the representation $\rho$ obtained by conjugation $\rho$ \cite{Weil1964Remarks}.  Namely, the coboundaries are the tangent vectors to derivatives of $g_{t}\rho g_{t}^{-1}$ where $g_{t}$ is some smooth path in $G$ based at the identity.  Thus, we may think of $H^{1}(\Gamma; \g_{\Ad \rho})$ as the tangent space to the \emph{character variety} at the equivalence class of representations of $\rho$ up to conjugation.\\
\\
Let us now specialize to the case where $\Gamma$ is finitely presented, as below.  
\begin{equation}\label{eq-defone}
\Gamma = \langle a_{1}, a_{2}, \hdots, a_{n} \, | \, r_{1},r_{2},\hdots, r_{n-1}\rangle 
\end{equation}
The machinery developed by Fox and his free-differential calculus aids us in evaluating crossed homomorphisms on any word $w(a_{1},a_{2},\hdots, a_{n})$ in terms of the crossed homomorphisms values on the generators.  Recall that a \emph{derivation} on $F_{n} \lra \Z F_{n}$ is a map satisfying $d(ab) = d(a) + ad(b)$ for all $a,b \in F_{n}$.  In \cite{Fox1953Free}, Fox showed the set of all such derivations of $F_{n}$ is in fact a \emph{free} $\Z F_{n}$-module generated by the partial derivative derivations $\del/\del a_{j}$ defined by 
\begin{equation}\label{eq-partial}
\frac{\del}{\del a_{j}} : F_{n} \lra \Z F_{n} \text{ via } \frac{\del a_{i}}{\del a_{j}}  = \delta_{i,j}
\end{equation}
Given a word $w \in F_{n}$, we denote the partial derivative of $w$ with respect to $a_{i}$ via $\frac{\del w}{\del a_{i}}$.  It is edifying to consider an example or two.  
\begin{example}\label{ex-partialinv}
The derivative of 1 is 0.
\begin{equation*}
\frac{\del}{\del a} 1 = \frac{\del}{\del a} 1\cdot 1 = \left(\frac{\del}{\del a} 1\right)+ 1 \left(\frac{\del}{\del a} 1\right) \implies \left(\frac{\del}{\del a} 1\right) = 0
\end{equation*}
The derivative of the $a^{-1}$ with respect to $a$ is $-a^{-1}$.   
\begin{equation*}
0 = \frac{\del}{\del a} a^{-1}a = \frac{\del}{\del a} a^{-1} + a^{-1}\frac{\del}{\del a} a \implies \frac{\del}{\del a} a^{-1} = -a^{-1} 
\end{equation*}
\end{example}
\begin{example}\label{ex-partialder}
Let $w = aba^{-1}b^{-1} \in F_{2}$.  Then
\begin{align*}
\frac{\del w}{\del a} &=1 + a \frac{\del}{\del a}\left( ba^{-1}b^{-1}\right) = 1 + a\left( \frac{\del b}{\del a} + b\frac{\del }{\del a} \left(a^{-1}b^{-1}\right) \right)\\
&=1 + ab\left(-a^{-1} + \frac{\del }{\del a}b^{-1}\right) = 1-aba^{-1}
\end{align*}
Similarly, 
\begin{equation*}
\frac{\del w}{\del b} = a-aba^{-1}b^{-1}
\end{equation*}
\end{example}
One can see from the above examples that the algebra of crossed homomorphisms and Fox-derivatives are intimately related.  In fact, given a finitely presentable group $\Gamma$, we have that for any word $w$ expressed in terms of the generators $a_{1},\hdots, a_{n}$, that
\begin{equation*}
f(w) = \sum_{i=1}^{n} \frac{\del w}{\del a_{i}} f(a_{i})
\end{equation*}
One can verify this simply by noting the right-hand side of the above equation is indeed a cocycle.  Moreover, it agrees with $f$ on the generators, and thus the two are equal.  This is the proof used in Proposition 3.2 as by Goldman \cite{Goldman2020Parallelism}.  We conclude this section by observing that in the case where $f$ is a coboundary, we have the following fundamental identity
\begin{equation*}
1-w = \sum_{i=1}^{n} \frac{\del w}{\del a_{i}} (1-a_{i})
\end{equation*}
\end{subsection}


\begin{subsection}{Twisted Alexander Polynomial}\label{subsec-twistedalex}
In this subsection we briefly recount and provide some examples of the \emph{twisted Alexander Polynomial}.  To the author's knowledge both Lin and Wada independently defined the twisted Alexander Polynomial for any finitely presentable group $\Gamma$ equipped with a representation into $\GL(l,D)$ where $D$ is a Euclidean domain accompanied by a surjection of $\Gamma$ onto a finite rank free-abelian group \cite{Lin2001RepresentationsOK} \cite{WADA1994241}.  This polynomial has largely been motivated and studied by knot complements but for our purposes we will be interested in the family of groups $\Gamma$ that admit the a presentation as in Equation \ref{eq-defone}.  Such groups are said to have \emph{deficiency one}.  \\
\\
Let $\rho: \Gamma \lra \GL(l,F)$ be a representation of $\Gamma$ into a vector space of dimension $l$ over $F$; typically for our purposes $F$ will either be $\R$ or $\C$.  Such a representation induces a homomorphism of group rings, which we abusively also write as $\rho: \Z\Gamma \lra \Z \GL(l,F) \subset \M(l,F)$ where $\M(l,F)$ denotes the ring of all $l\times l$ matrices with entries in $F$.  \\
\\
We assume further we are provided a \emph{surjective} group homomorphism $\alpha: \Gamma \lra \langle x \rangle$ where $\langle x \rangle $ is an infinite cyclic group isomorphic to $\Z$.  Note the group ring of $\langle x \rangle$ is readily seen to be isomorphic to the Laurent polynomial ring $\Z[x,x^{-1}]$.  The group homomorphism $\alpha$ induces a ring homomorphism which we denote $\alpha: \Z\Gamma \lra \Z[x,x^{-1}]$.  We take the tensor product of $\rho$ and $\alpha$ to form the ring homomorphism
\begin{equation}\label{eq-tensormap}
\rho \otimes \alpha : \Z\Gamma \lra \M(l,F)\otimes_{\Z} \Z[x,x^{-1}] \simeq M\left(l, F[x,x^{-1}]\right)
\end{equation}
Finally, we take the quotient map induced by the presentation as in Equation \ref{eq-defone}, $F_{n} \lra \Gamma$, which induces a map on the level of group rings $\Z F_{n} \lra \Z \Gamma$, and (post)-compose this group ring map with the above map in Equation \ref{eq-tensormap}.  This yields a ring homomorphism which we denote by
\begin{equation}\label{eq-phi}
\Phi: \Z F_{n} \lra M\left(l, F[x,x^{-1}]\right)
\end{equation}
We now introduce a $(n-1)\times n$-matrix $A_{\rho,\alpha}$ which we call the \emph{Twisted Alexander Matrix}.  
\begin{definition}\label{def-twistedmat}
Let $\mathcal{A}_{\rho,\alpha}$ denote the $(n-1)\times n$-matrix whose $(i,j)$-th entry is an $l \times l$-matrix 
\begin{equation*}
A_{i,j} = \Phi\left(\frac{\del r_{i}}{\del a_{j}}\right) \in M\left(l, F[x,x^{-1}]\right)
\end{equation*}
where $\del r_{i}/\del a_{j}$ is the Fox-derivative of the relation $r_{i}$ as in Equation \ref{eq-defone}.  
\end{definition}
We denote by $\mathcal{A}_{\rho,\alpha,k}$ the $(n-1)\times(n-1)$-matrix obtained from $\mathcal{A}_{\rho,\alpha}$ by removing the $k$-th column.  This is a square $(n-1)\times(n-1)$-matrix whose \emph{entries} are $l\times l$-matrices with coefficients in $F[x,x^{-1}]$, which we may consider a $l(n-1)\times l(n-1)$-matrix with entries in $F[x,x^{-1}]$.  \\
\\
Observe that there exists a $k$ such that $\alpha(a_{k}) \neq 0 \in \Z$ as $\alpha$ is surjective onto $\Z$.  Then $\Phi(1-a_{k}) = 1-\alpha(a_{k})\rho(a_{k})$ is some non-trivial Laurent polynomial in $x$, thus $\det \Phi(1-a_{k}) \neq 0$ for some $k$.  With this we are able to define the twisted Alexander polynomial.  
\begin{definition}\label{def-twistedalex}
Let $\rho: \Gamma \lra \GL(l,F)$ be a representation of a finitely presented group with deficiency one, $\Gamma = \langle a_{1}, a_{2}, \hdots, a_{n} \, | \, r_{1},r_{2},\hdots, r_{n-1}\rangle$, into an $l$-dimensional vector space over the field $F$.  Let $\alpha: \Gamma \lra \Z$ be a surjective homomorphism.  Let $\Phi: F_{n} \lra M(l,F[x,x^{-1}])$ be the composition of the group ring maps $\Z F_{n} \lra Z\Gamma$ with $\rho \otimes \alpha: \Z \Gamma \lra M(l,F)\otimes_{\Z}\Z[x,x^{-1}] \simeq M(l,F[x,x^{-1}])$.  Form the twisted Alexander matrix $\mathcal{A}_{\rho,\alpha}$ as in Definition \ref{def-twistedmat}, and choose a $k$ such that $\det \Phi(1- a_{k}) \neq 0$.  We define the twisted Alexander polynomial to be 
\begin{equation}
\Delta_{\rho,\alpha}(x) = \frac{\det \mathcal{A}_{\rho,\alpha, k}}{\det \Phi\left(1-a_{k}\right)}
\end{equation}
where $\mathcal{A}_{\rho,\alpha,k}$ is the twisted Alexander matrix with its $k$-th column removed.
\end{definition}
We provide several remarks before providing examples.  The twisted Alexander polynomial is only well-defined up to multiples of elements in $F$ and powers of $x^{ls}$ for integers $s \in \Z$; however one may rationalize to obtain an honest to goodness polynomial.  The twisted Alexander polynomial also does not depend on the choice of $k$ for which $\det \Phi(1- a_{k}) \neq 0$ \cite{WADA1994241}.    \\
\\
To the author's knowledge the twisted Alexander polynomial has been most widely studied in the context of hyperbolic geometry.  Specifically, given a hyperbolic knot complement $M = S^{3}\setminus K$, Thurston's Hyperbolization Theorem says that $M$ admits a unique complete hyperbolic structure \cite{Thurston2022Geometry}.  Such a structure is determined by a discrete and faithful representation of $\pi_{1}(M) \lra \PSL(2,\C)$, which one may lift to $\SL(2,\C)$ by choice of a spin structure \cite{CULLER198664}.  A Mayer-Vietrois sequence argument readily shows that $H_{1}(M)$ is isomorphic to $\Z$, and thus the abelinization map of $\pi_{1}(M)$ provides our surjective homomorphism $\alpha$.  Dunfield, Friedl, and Jackson infer topological data from the polynomial in this context \cite{Dunfield2012Twisted} whereas Boden and Friedl \cite{BODEN_FRIEDL_2014} illustrate the mysteriousness of these polynomials by illustrating their tensor products do not behave in an expected manner.  
\\
\\
With all the necessary machinery established, we move on to some examples regarding the trefoil knot with presentation given by $\Gamma = \langle a,b \, | \, a^{2}b^{-3} \rangle$.  
\begin{example}\label{trivtrefoil}
We first consider the trivial representation of the trefoil knot complement group $\rho: \Gamma \lra \GL(1,\R)$ and we take as an abelinization map $\alpha: \Gamma \lra \langle x \rangle$, via $\alpha(a) = x^{3}$ and $\alpha(b) = x^{2}$.  We need to calculate the Fox-derivatives of the relation $r_{1} = a^{2}b^{-3}$.  They are given below.
\begin{equation}\label{eq-foxdertref}
\frac{\del r_{1}}{da} = 1 + a \phantom{ and } \frac{\del r_{1}}{db} = -a^{2}b^{-1} - a^{2}b^{-2}-1
\end{equation}
Our twisted Alexander Matrix is given by
\begin{equation*}\label{eq-trivalex}
\mathcal{A}_{\rho,\alpha}(x) = \left(
\begin{array}{cc}
1+x^{3} & -x^{4}-x^{2}-1
\end{array}
\right)
\end{equation*}
If we take $\mathcal{A}_{\alpha,\rho,2}$ and $\Phi(1-a) = 1-x^{3}$, our twisted Alexander polynomial is given by
\begin{equation*}
\Delta_{\rho,\alpha} = \frac{\det \mathcal{A}_{\rho,\alpha,2}}{\det \Phi(1-a)} =  \frac{1-x-x^{2}}{x-1}
\end{equation*}
which is the typical \emph{untwisted} Alexander polynomial.  One can readily verify that had we chosen the first column of the twisted Alexander matrix, we would end up with the same result, up to sign.  
\end{example}
\begin{example}\label{heusenertrefoil}
This next example follows up on the work of Heusener where he classifies all irreducible representations of $\Gamma = \langle a,b \, | \, a^{2}b^{-3} \rangle$ into $\SL(3,\R)$ \cite{gomtrie2016SomeRR}.  Let $\omega$ denote a primitive 3rd root of unity and consider a representation $\rho_{s,t}: \Gamma \lra \SL(3,\R)$ via
\begin{equation}
\rho(a) = \left(
\begin{array}{ccc}
1 & 0 & 0 \\
s & -1 & 0 \\
t & 0 & -1
\end{array}
\right) \text{ and }
\rho(b) = \left(
\begin{array}{ccc}
1 & \omega-1 & \omega^{2}-1 \\
0 & \omega & 0 \\
0 & 0 & \omega^{2}
\end{array}
\right)
\end{equation}
where $s,t \in \C$.  We use the same abelinization map as in Example \ref{trivtrefoil}, and calculate the twisted Alexander Matrix.
\begin{scriptsize}
\begin{equation*}\label{eq-heusalex}
\mathcal{A}_{\rho_{s,t},\alpha} = \left(
\begin{array}{cc}
	\left(
		\begin{array}{ccc}
		 x^3+1 & 0 & 0 \\
		 s x^3 & 1-x^3 & 0 \\
		 t x^3 & 0 & 1-x^3 \\
		\end{array}
	\right)
	&
	\left(
		\begin{array}{ccc}
		 -1-x^{2}-x^{4} & (1-\omega)x^{2} + (1-\omega^{2})x^{4}& (1-\omega^{2})x^{2} + (1-\omega)x^{4} \\
		 0 & -1-\omega x^{2} - \omega^{2}x^{4} & 0 \\
		 0 & 0 & -1-\omega^{2}x^{2}-\omega x^{4} \\
		\end{array}
	\right)
\end{array}
\right)
\end{equation*}
\end{scriptsize}
Taking $\mathcal{A}_{\rho_{s,t},\alpha,1} $ and 
\begin{equation*}
\Phi(1-b) =  \left(
\begin{array}{ccc}
1-x^{2} &(1-\omega)x^{2} & (1-\omega^{2})x^{2} \\
0 & 1-\omega x^{2} & 0 \\
0 & 0 & 1-\omega^{2}x^{2}
\end{array}
\right)
\end{equation*}
we see that 
\begin{equation*}
\Delta_{\rho_{s,t},\alpha}(x) = \frac{\det \mathcal{A}_{\rho,\alpha,2}}{\det \Phi(1-a)} = \frac{(1+x^{3})(1-x^{3})^{2}}{(1-x^{2})(1-\omega x^{2})(1-\omega^{2}x^{2})} = \frac{(1+x^{3})(1-x^{3})^{2}}{(1-x^{6})} = 1-x^{3}
\end{equation*}
Note that while the original representation has two whole complex parameters of freedom, the twisted Alexander polynomial is independent of them.  This particular example illustrates how coarse the invariant can be in its ability to distinguish certain representations.  For generic $(s,t) \in \C^{2}$, $\rho_{s,t}$ is irreducible, however there are three affine lines for which $\rho_{s,t}$ is reducible \cite{gomtrie2016SomeRR}.  Moreover, at the three points of intersections of these lines, $\rho_{s,t}$ fixes a complete flag and has the character of a diagonal representation.  However, none of these distinguishing factors is detected by the twisted Alexander polynomial.  
\end{example}

\end{subsection}

\begin{subsection}{Infinitesimal Projective Rigidity}\label{ssipr}
In this section we review some concepts regarding the notion of infinitesimal projective rigidity.  To motivate the conversation, let $M$ be a closed orientable 3-manifold with fundamental group $\Gamma$.  If $M$ admits a hyperbolic structure, then there is a discrete and faithful representation, $\rho: \Gamma \lra \PSL(2,\C)$ of its fundamental group into the orientation preserving isometries of hyperbolic 3-space.  Mostow's famous rigidity theorem states that any two such representations are in fact conjugate in $\PSL(2,\C)$; i.e. given another discrete faithful representation $\sigma: \Gamma \lra \PSL(2,\C)$, there exists a $g \in \PSL(2,\C)$ such that $\rho = g\sigma g^{-1}$.  \\
\\
In these representations we are identifying hyperbolic 3-space, $H^{3}$, with upper-half space, and its orientation preserving isometry group with $\PSL(2,\C)$.  One could instead use the \emph{Klein} model of $H^{3}$.  Let $\R^{3,1}$ denote $\R^{4}$ equipped with the non-degenerate symmetric bilinear form $(\phantom{x}, \phantom{y})$ defined below
\begin{equation}\label{eq-lorentz}
(\mathbf{x}, \mathbf{y}) = x_{1}y_{1} + x_{2}y_{2} + x_{3}y_{3} - x_{4}y_{4} \text{ for all } \mathbf{x},\mathbf{y} \in \R^{4}
\end{equation}
This is the Lortenzian metric of signature $(3,1)$.  If we consider the hyperboloid $H = \{\mathbf{x} \in \R^{3,1} \, | \, (\mathbf{x}, \mathbf{x}) = -1 \}$ and projectivize it under $\R^{4}\setminus 0 \lra \RP^{3}$, its image is the Klein model for hyperbolic 3-space which we denote by $K^{3}$.  Its group of isometries is given by $\PSO(3,1)$, the projectivization of $\SO^{0}(3,1)$.  \\
\\
Correspondingly, another means of obtaining a hyperbolic 3-manifold is by an atlas of charts into $K$ whose coordinate changes are restrictions of isometries of $\PSO(3,1)$.  Thus, the holonomy representation of $\Gamma$ may instead be taken into $\PSO(3,1)$.  As such, we see every hyperbolic 3-manifold (with or without boundary) inherits a natural \emph{projective structure}. \\
\\
What is less well known is whether the projective structure of $M$ deforms in the larger group of projectivities $\PGL(4,\R)$.  That is to say, does there exist a smooth path of representations $\rho_{t}: \Gamma \lra \PGL(4,\R)$ such that $\rho_{0} = \rho$ and $\rho_{t}$ is not simply obtained by conjugation?  Heuristically, the work of several authors, Cooper-Long-Thistlethwaite, Daly, Heusener-Porti, Ballas-Danciger-Lee, says this behavior is quite rare in the $\emph{closed}$ case.  In Cooper-Long-Thistlethwaite they conducted a large scale census of approximately $4,500$ manifolds and determined that at most 61 or so of them admit such deformations \cite{Cooper2006Computing}.  In Heusener-Porti, they showed that the closed hyperbolic 3-manifolds obtained via Dehn-filling of the figure-eight knot complement for sufficiently large $(p,q)$ do not admit such deformations \cite{Heusener2011Infinitesimal}.  In Daly, we showed that the Dehn-fillings of the figure-eight knot complement for many small $(p,q)$ similarly do not admit these deformations \cite{Daly}.  \\
\\
Many of the arguments rely heavily the notion of \emph{infinitesimal projective rigidity} which we recall now.  Let $M$ denote an orientable \emph{closed} hyperbolic 3-manifold with fundamental group $\Gamma$ and holonomy representation $\Gamma \lra \PSO(3,1) \subset \PGL(4,\R)$.  The Lie-algebra of $\PGL(4,\R)$ may be identified with $\slf(4,\R)$ in the standard manner.  We endow $\slf(4,\R)$ with the structure of a $\Gamma$-module by post-composing the holonomy representation $\rho: \Gamma \lra \PGL(4,\R)$ with the adjoint representation $\Ad: \PGL(4,\R) \lra \Aut (\slf(4,\R))$.  Denote the corresponding $\Gamma$-module by $\slf(4,\R)_{\Ad \rho}$. 
\begin{definition}\label{ifr}
We say that such a manifold $M$ is infinitesimally projectively rigid if and only if $H^{1}(\Gamma; \slf(4,\R)_{\Ad \rho}) = 0$.
\end{definition}
As mentioned in Section \ref{subsec-twistedcohom} one may think of $H^{1}(\Gamma; \slf(4,\R)_{\Ad \rho})$ as the tangent space to the character variety of $\Gamma$ into $\PGL(4,\R)$ at the point corresponding to holonomy representation $\rho$.  Recall the Ehresmann-Weil-Thurston Principle which loosely states that equivalent $(G,X)$-structures on a manifold $M$ are in local correspondence with the space of representations of the fundamental group of $M$ into $G$ modulo conjugation.  So a deformation of the hyperbolic structure of $M$ to a non-equivalent projective structure will induce non-conjugate holonomy representations.  If this deformation is done in a smooth manner, we will obtain a smooth path of holonomy representations $\rho_{t}$ starting at the hyperbolic structure $\rho_{0} = \rho$.  This will correspond to a smooth path in the representation variety.  Our interest though is in \emph{non-equivalent} projective structures, so we do not want to consider projective structures on $M$ that are equivalent to the original hyperbolic structure.  Equivalent structures will give rise to conjugate holonomy representations $\rho_{t}$ and for sufficiently close ones to the hyperbolic structure, we may find a smooth path $g_{t}$ of projectivities such that $\rho_{t} = g_{t}\rho g_{t}^{-1}$ where $\rho$ is the original holonomy representation.\\
\\
We say that the manifold $M$ is \emph{locally projectively rigid} if every such smooth path of representations $\rho_{t}$ starting at $\rho_{0} = \rho$ is trivial, in the sense that it arises from conjugation as $\rho_{t} = g_{t}\rho g_{t}^{-1}$.  This is to say, there are no non-trivial deformations of the hyperbolic structure in question.  Weil's Infinitesimal Rigidity Theorem states that infinitesimal rigidity implies local rigidity.  Thus if this group is zero, then all deformations are trivial and there is no non-trivial deformation of the hyperbolic structure in question.  \\
\\
A related notion is that of \emph{infinitesimal projective rigidity rel cusp(s)}.  In this context $M$ is a \emph{cusped} orientable hyperbolic 3-manifold, thus its boundary $\del M$ is a finite union of tori.  Let $\rho$ denote the holonomy representation of the \emph{complete} hyperbolic structure of $M$ which is unique up to conjugation by Mostow Rigidity and $\Gamma$ denote its fundamental group.  Due to the fact that $M$ is a quotient of hyperbolic 3-space, which is contractible, we may identify the abstract group theoretic cohomology of the $\Gamma$-module, $\slf(4,\R)_{\Ad \rho}$, with the twisted \emph{topological} cohomology associated to the representation of $\Gamma$ given by $\Ad \rho$.  See Hatcher for more details \cite{Hatcher2002Algebraic}.  
\begin{definition}\label{defiprrc}
An orientable cusped hyperbolic 3-manifold $M$ is infinitesimally projectively rigid rel cusp(s) if and only if $i^{*}: H^{1}(M; \slf(4,\R)_{\Ad \rho}) \lra H^{1}(\del M; \slf(4,\R)_{\Ad \rho}) $ induced by $\del M \subset M$ is an injective homomorphism where $i^{*}$ is the map on cohomology induced by inclusion of $\del M \subset M$.  
\end{definition}
The strength of this hypothesis has been used with great reward by many authors including Heusener-Porti, Ballas-Danciger-Lee, and Ballas-Danciger-Lee-Marquis \cite{Heusener2011Infinitesimal}, \cite{Ballas2018Convex}, \cite{BDLM}.  For example, the first pair of authors use this hypothesis along with other cohomological techniques to prove their result regarding infinitesimal projective rigidity of infinitely many Dehn-surgeries on the figure-eight knot complement.  The second trio of authors use the hypothesis to show that manifolds that satisfy Definition \ref{defiprrc} have properly convex projective structures on their \emph{doubles}.  In addition they showed that the corresponding point in the character variety associated to the complete hyperbolic structure is \emph{smooth}.  \\
\\
In Heusener-Porti and Hodgson they use the so-called `Half-Lives Half-Dies' theorem which states that the map $i^{*}: H^{1}(M; \slf(4,\R)_{\Ad \rho}) \lra H^{1}(\del M; \slf(4,\R)_{\Ad \rho})$ above in Definition \ref{defiprrc} has image equal to $\frac{1}{2} \dim H^{1}(\del M; \slf(4,\R)_{\Ad \rho}) = 3$ \cite{Heusener2011Infinitesimal} \cite{Hodgson}.  This dimension equality follows from a straight-forward application of Poincar\'{e} duality and an Euler characteristic argument.  By Poincar\'{e} duality, one has $\dim H^{0}(\del M; \slf(4,\R)_{\Ad \rho}) = \dim H^{2}(\del M; \slf(4,\R)_{\Ad \rho})$, the former of which is equal to 3 by determining the dimension of the infinitesimal centralizer of the fundamental group of the boundary torus in $\SL(4,\R)$.  Because the Euler characteristic of the torus vanishes, $\dim H^{1}(\del M; \slf(4,\R)_{\Ad \rho})$ $= 2\dim H^{0}(\del M; \slf(4,\R)_{\Ad \rho})$.  Thus the dimension of the image of $i^{*}$ inside $H^{1}(\del M; \slf(4,\R)_{\Ad \rho})$ is equal to 3.  One can see Heusener-Porti, Ballas-Danciger-Lee, or Daly for more details \cite{Heusener2011Infinitesimal} \cite{Ballas2014Deformations} \cite{Daly}.  Summarizing these calculations and the definition of infinitesimal projective rigidity rel cusp, we have the following lemma.
\begin{lemma}\label{lem1}
Let $M$ be a hyperbolic once-punctured torus bundle with holonomy $\rho: \Gamma \lra \PSO(3,1) \subset \PGL(4,\R)$.  Then $M$ is infinitesimally projectively rigid rel cusp if and only if $\dim H^{1}(M; \slf(4,\R)_{\Ad \rho}) = 3$.  
\end{lemma}
The current convention regarding verifying whether or not a cusped orientable hyperbolic 3-manifold is infinitesimally projectively rigid rel cusp(s) is direct calculation.  In Heusener-Porti they prove that all but finitely many punctured torus bundles with tunnel number one satisfy this condition, however their result does not provide a specific bound.  To the author's knowledge there is not an explicit family known of cusped orientable hyperbolic 3-manifolds that satisfy the condition of infinitesimal projective rigidity rel cusp(s).  It is the goal of this paper to provide equivalent, and sufficient, conditions to verify infinitesimal projective rigidity rel cusp in the case of once-punctured torus bundles.   
\end{subsection}

\end{section}

\begin{section}{Cohomology Calculations}\label{sectioncc}
In this section we aim to provide cohomology calculations regarding the twisted cohomology of the hyperbolic structure on $M$.  We begin with some group theoretic and geometric preliminaries.  
\begin{subsection}{Once-Punctured Torus Bundles}\label{ssoptb}
For our purposes, we will be interested in hyperbolic once-punctured torus bundles which are fiber bundles $M$ of the form $F \lra M \lra S^{1}$ where $F$ is a once-punctured torus.  The fundamental group of $F$ is isomorphic to a free group of rank two, which we denote by $F_{2}$ and generated by elements $a$ and $b$.  The base circle's fundamental group is generated by an element $x$.  The monodromy of the action of $x$ induces an outer automorphism of $\pi_{1}(F)$ which entirely determines the once-punctured torus bundle up to conjugacy, thus the topological structure of $M$ is determined by a single element $\phi \in \SL(2,\Z)$ satisfying $|\tr \phi|> 2$.  Moreover, by a theorem of Murasugi, $\phi$ admits a decomposition as
\begin{equation*}
\phi = i^{\epsilon}R^{n_{1}}L^{n_{2}}\hdots L^{n_{2k}}
\end{equation*}
where
\begin{equation*}
i = \left(
\begin{array}{cc}
-1 & 0 \\
0 & -1
\end{array}
\right) \phantom{==}
R = \left(
\begin{array}{cc}
1 & 1 \\
0 & 1
\end{array}
\right) \phantom{==}
L = \left(
\begin{array}{cc}
1 & 0 \\
1 & 1
\end{array}
\right)
\end{equation*}
The group theoretic analogue of the fiber bundle structure is the split short exact sequence of groups in Equation \ref{eqsesgroups}.  
\begin{equation}\label{eqsesgroups}
F_{2} \lra \pi_{1}(M) \lra F_{1}
\end{equation}
Here $F_{k}$ denotes the free group on $k$-generators.  The group $F_{2}$ is a normal subgroup corresponding to the fundamental group of the once-punctured torus fiber $F$ and the $F_{1}$ quotient is the fundamental group of the base circle.  For a fixed hyperbolic once-punctured torus bundle $M$ with monodromy $\phi$, we obtain the presentation of the fundamental group of $M$ as in Equation \ref{eqpres}
\begin{equation}\label{eqpres}
\pi_{1}(M) = \langle a,b,x \, | \, xax^{-1} = \phi(a) \text{ and } xbx^{-1} = \phi(b) \rangle
\end{equation}
The subgroup $a,b$ in Equation \ref{eqpres} generate a free group of rank 2 and corresponds to the fundamental group of the fiber $F$ whereas $x$ is the image of a section of the base circle and represents a meridian of the toroidal boundary.  A longitude is given by the boundary of the once-puncture fiber and expressed in terms of our generators as the commutator of $a$ and $b$ which we denote by $l = aba^{-1}b^{-1}$.  \\
\\
Both J{\o}rgensen and Thurston showed that a hyperbolic once-punctured torus bundle admits a unique complete hyperbolic structure, and consequently, a discrete and faithful representation of its fundamental group $\Gamma$ into $\PSL(2,\C)$ \cite{Jorgensen} \cite{Thurston2022Geometry}.  This is the \emph{holonomy representation} of the hyperbolic structure on $M$ which we denote by $\rho: \Gamma \lra \PSL(2,\C)$.  In theory, one can calculate this representation explicitly.  To explain how, we introduce a bit of notation.  Let $A = \tr\rho(a), B = \tr\rho(b)$, $C = \tr \rho(ab)$ which by Equation \ref{eqpres} satisfy the trace equations in Equation \ref{eqtraceeqs}.
\begin{equation}\label{eqtraceeqs}
A = \tr \rho( \phi(a)), \phantom{==} B = \tr \rho(\phi(b)), \phantom{==} A^{2} + B^{2} + C^{2} = ABC
\end{equation}
The last equation is a consequence of the completeness of $M$ and comes from the fact that the holonomy about the longitude $l = aba^{-1}b^{-1}$ is parabolic.  In fact if we begin with a triple $A,B,C \in \C$ which satisfies Equation \ref{eqtraceeqs}, then we may furnish a representation of $\Gamma$ by first considering the discrete and faithful representation of $F_{2}$ given by the Maclachlan-Reid representation of the free group on 2-generators (Section 4.6 \cite{Maclachlan2003Arithmetic}) in Equation \ref{eqfree2} below.
 \begin{equation}\label{eqfree2}
 \rho(a) = \frac{1}{C}\left(
 \begin{array}{cc}
 AC-B & A/C \\
 AC & B
 \end{array}
 \right),\phantom{==}
  \rho(b) = \frac{1}{C}\left(
 \begin{array}{cc}
 BC-A & -B/C \\
 -BC & A
 \end{array}
 \right)
 \end{equation}
Using the fact that $\rho(x)$ is parabolic and the relations in Equation \ref{eqpres}, we may then solve for $\rho(x)$ to obtain the geometric representation of $\Gamma$ into $\PSL(2,\C)$.  Solving the system of equations in Equation \ref{eqtraceeqs} is typically a difficult problem, however there are low-complexity scenarios where we may explicitly write out the roots of these equations.  This is typically done by employing the so-called \emph{trace relations} originally studied by Fricke and Klein and later by Goldman \cite{1914Vorlesungen}, \cite{Goldman2009Trace}.  We provide an example below.
\begin{example}\label{exLLRR}
Let $M$ be the hyperbolic once-punctured torus bundle with monodromy $\phi = L^{2}R^{2}$ as below.
\begin{equation}
L^{2}R^{2} = \left(
\begin{array}{cc}
1 & 0 \\
1 & 1
\end{array}
\right)^{2}
 \left(
\begin{array}{cc}
1 & 1 \\
0 & 1
\end{array}
\right)^{2} = 
 \left(
\begin{array}{cc}
1 & 2 \\
2 & 5
\end{array}
\right)
\end{equation}
Let $F_{2}$ be the fundamental group of the once-puncture torus generated by $a$ and $b$ and $x$ be the meridian which acts by $\phi$ on $a$ and $b$ when conjugated.  Recall that $L$ and $R$ act on $F_{2}$ via the automorphisms below.
\begin{align*}
L(a) &= ab \text{ and } L(b) = b\\
R(a) &= a \text{ and } R(b) = ba
\end{align*} 
Thus $\Gamma$ admits a presentation of the form below.
\begin{equation}\label{eqtracerels}
\Gamma = \langle a,b,x \, | \, xax^{-1} = ab^{2} \text{ and } xbx^{-1} = bab^{2}ab^{2} \rangle
\end{equation}
Let $A,B,C$ denote the traces as defined in the previous paragraphs.  To solve the equations as defined by Equation \ref{eqtraceeqs}, we implement to trace-relations as in Maclachlan and Reid \cite{Maclachlan2003Arithmetic}.  Specifically for any matrices $a,b \in \SL(2,\C)$ we have the below equalities.  
\begin{align*}
\tr ab  = \tr a\, \tr b - \tr ab^{-1} \phantom{==} \tr a^{2} = \tr^{2} a - 2 \phantom{==} \tr a= \tr a^{-1}
\end{align*}
So for example we may reduce 
\begin{equation}
\tr ab^{2}  = \tr ab\tr b - \tr a = CB-A
\end{equation}
Similarly we have
\begin{align*}
\tr (bab^{2})(ab^{2}) &= \tr bab^{2}\tr ab^{2} - \tr b \\
&= \left(\tr ab^{2} \tr b -\tr ab\right) \tr ab^{2} - \tr b \\
&= \left((CB-A)B - C\right)(CB-A) - B
\end{align*}
Thus we have the polynomial equations below.
\begin{equation}
(CB-A) = A \phantom{==} \left((CB-A)B - C\right)(CB-A) - B = B \phantom{==} A^{2} + B^{2} + C^{2} = ABC
\end{equation}
Solving for $B$ yields that $B = -\sqrt{2-2i}$ and $A$ and $C$ are given by
\begin{equation*}
A = -\frac{\sqrt{2}B}{\sqrt{B^{2}-2}} \phantom{==} C = -\frac{2\sqrt{2}}{\sqrt{B^{2}-2}}
\end{equation*}
These calculations are carried out explicitly in \texttt{L2R2}\underscore\texttt{example.nb} on the author's webpage.  Equation \ref{eqfree2} yields the holonomy of our once-punctured torus as
\begin{footnotesize}
\begin{equation}\label{eqf2holL2R2}
\rho(a) = \left(
\begin{array}{cc}
 -\frac{\left(\frac{1}{2}-\frac{i}{2}\right) \left(\sqrt{2-2 i}-\frac{4
   \sqrt{2}}{(1-i)^{3/2}}\right)}{\sqrt{2}} & \frac{1}{4} (1-i)^{3/2} \\
 \frac{2}{\sqrt{1-i}} & \frac{1}{2} (1-i)^{3/2} \\
\end{array}
\right) \phantom{==}
\rho(b) = \left(
\begin{array}{cc}
 -\sqrt{\frac{1}{2}-\frac{i}{2}} & \frac{(1-i)^{5/2}}{4 \sqrt{2}} \\
 \sqrt{2-2 i} & -\sqrt{\frac{1}{2}-\frac{i}{2}} \\
\end{array}
\right)
\end{equation}
\end{footnotesize}
Solving for $\rho(x)$ is a matter of solving the relations imposed by $xax^{-1} = ab^{2} \text{ and } xbx^{-1} = bab^{2}ab^{2}$ and that $\det \rho(x) = 1$.  Doing so yields that
\begin{equation*}
\rho(x) = \left(
\begin{array}{cc}
 -1 & i \\
 0 & -1 \\
\end{array}
\right)
\end{equation*}
Up to conjugation, this is the unique complete structure on the hyperbolic once-punctured torus bundle $M$.  With the explicit holonomy representation, we may use the isomorphism between $\PSL(2,\C)$ and $\PSO(3,1)$ as described in Cooper-Long-Thistlethwaite to yield the holonomy representation into $\PSO(3,1)$ given below \cite{Cooper2006Computing}.
\begin{scriptsize}
\begin{equation*}
\rho(a) = \left(
\begin{array}{cccc}
 -\frac{1}{\sqrt{2}} & -\frac{1}{\sqrt{2}} & \sqrt{2} & -\sqrt{2} \\
 \frac{3}{\sqrt{2}} & -\frac{1}{\sqrt{2}} & -\frac{7}{2 \sqrt{2}} & \frac{9}{2 \sqrt{2}}
   \\
 \frac{1}{\sqrt{2}} & \frac{1}{2 \sqrt{2}} & \frac{7}{8 \sqrt{2}} & -\frac{1}{8 \sqrt{2}}
   \\
 \frac{3}{\sqrt{2}} & -\frac{1}{2 \sqrt{2}} & -\frac{31}{8 \sqrt{2}} & \frac{41}{8
   \sqrt{2}} \\
\end{array}
\right) 
\rho(b) = \left(
\begin{array}{cccc}
 \frac{1}{\sqrt{2}} & \frac{1}{\sqrt{2}} & -\frac{1}{2 \sqrt{2}} & -\frac{1}{2 \sqrt{2}}
   \\
 \frac{1}{\sqrt{2}} & \frac{1}{\sqrt{2}} & \sqrt{2} & -\sqrt{2} \\
 \frac{1}{2 \sqrt{2}} & -\sqrt{2} & -\frac{9}{8 \sqrt{2}} & \frac{15}{8 \sqrt{2}} \\
 -\frac{1}{2 \sqrt{2}} & -\sqrt{2} & -\frac{15}{8 \sqrt{2}} & \frac{25}{8 \sqrt{2}} \\
\end{array}
\right)
\rho(x) = \left(
\begin{array}{cccc}
 1 & 0 & 1 & 1 \\
 0 & 1 & 0 & 0 \\
 -1 & 0 & \frac{1}{2} & -\frac{1}{2} \\
 1 & 0 & \frac{1}{2} & \frac{3}{2} \\
\end{array}
\right)
\end{equation*}
\end{scriptsize}  
\end{example}
\end{subsection}

\begin{subsection}{Cohomology of Once-Punctured Torus Bundles}\label{sscohomology}
In this section we provide some calculations on once-punctured torus bundles and rephrase the condition of being infinitesimally projectively rigid rel cusp in terms of an action on the character variety of $F_{2}$.  We let $M$ denote a hyperbolic once-punctured torus bundle with monodromy $\phi \in \SL(2,\Z)$ and fundamental group $\Gamma$.  Let $\rho: \Gamma \lra \PSO(3,1) \subset \PGL(4,\R)$ denote the holonomy representation of $M$ associated to its unique complete hyperbolic structure.  Compose this representation with the adjoint representation of $\PGL(4,\R)$ to obtain the representation $\Ad \rho: \Gamma \lra \Aut(\slf(4,\R))$ as we did in Section \ref{ssipr} so that $\slf(4,\R)$ becomes a $\Gamma$-module.  The main object of our interest is $H^{1}(\Gamma; \slf(4,\R)_{\Ad \rho})$.  Because our manifolds $M$ are hyperbolic, their universal covers are contractible, and thus we enjoy an isomorphism between the twisted group cohomology of the $\Gamma$-module $\slf(4,\R)_{\Ad \rho}$ and the twisted \emph{topological} cohomology $H^{1}(M; \slf(4,\R)_{\Ad \rho})$ where the coefficient module is twisted by the representation $\Ad \rho: \Gamma \lra \Aut(\slf(4,\R)_{\Ad \rho})$ as stated in Section \ref{ssipr}.  Because of this isomorphism, we enjoy many of the tools from ordinary topological cohomology theory such as the Mayer-Vietoris sequence, excision, long-exact sequences of good pairs, etc.    \\
\\
As done in Johnson and Millson, let us split the coefficient module $\slf(4,\R) = \pso(3,1)\oplus \vv$ where $\vv$ is the 9-dimensional orthogonal complement of $\pso(3,1)$ under the Killing form of $\slf(4,\R)$ \cite{Johnson1987Deformation}.  We emphasize this splitting is as $\Gamma$-modules and not Lie-algebras.  Thus we have that 
\begin{equation}\label{eqhsplit}
H^{1}(\Gamma; \slf(4,\R)_{\Ad \rho}) = H^{1}(\Gamma; \pso(3,1)_{\Ad \rho})\oplus H^{1}(\Gamma; \vv_{\Ad \rho})
\end{equation}
The first piece represents deformations of the geometric representation of $\rho$ inside $\PSO(3,1)$.  This is well understood due to Thurston's Hyperbolic Dehn-Filling Theorem \cite{Thurston2022Geometry}.  In particular, we know that  $\dim H^{1}(\Gamma; \pso(3,1)_{\Ad \rho}) = 2$, see Kapovich for details \cite{Kapovich2010Hyperbolic}.  We are therefore largely interested in the latter piece, $H^{1}(\Gamma; \vv_{\Ad \rho})$.  Because of Lemma \ref{lem1}, the fact that $\dim H^{1}(\Gamma; \pso(3,1)_{\Ad \rho}) = 2$, and the splitting in Equation \ref{eqhsplit}, we see that infinitesimal projectivitiy rigidity rel cusp in our context is equivalent to $\dim H^{1}(\Gamma; \vv_{\Ad \rho}) = 1$.  This is a specific application of Corollary 5.4 in Heusener-Porti \cite{Heusener2011Infinitesimal}.     \\
\\
In what follows, to ease notation we suppress the subscript $\Ad \rho$ when the context is clear.  Because $M$ admits the structure of a fiber bundle, the base circle acts on the cohomology groups of the fiber, $H^{*}(F_{2}; \slf(4,\R))$.  As the wedge of two circles is a $K(F_{2},1)$, by the equivalence of topological and group cohomology for spaces with contractible universal covers, the cohomology groups $H^{*}(F_{2}; \slf(4,\R))$ vanish in degrees larger than 1.  Recall that $H^{0}(G;V)$ is isomorphic to the invariants of the group action of $G$ on $V$ \cite{Brown1982Cohomology}.  However, the holonomy representation of $\rho: F_{2} \lra \PSO(3,1) \subset \PGL(4,\R)$ is irreducible, thus the centralizer of $\rho(F_{2})$ inside $\PGL(4,\R)$ is discrete, and its infinitesimal centralizer is trivial so $H^{0}(F_{2}; \slf(4,\R)) = 0$.  These calculations lend themselves nicely to an application of the Lyndon–Hochschild–Serre spectral sequence which states that spectral sequence in Equation \ref{eqlhs} converges to the cohomology of $\Gamma$.  
\begin{equation}\label{eqlhs}
E_{2}^{p,q} = H^{p}(\Gamma/F_{2};H^{q}(F_{2};\slf(4,\R))) \implies H^{p+q}(\Gamma; \slf(4,\R))
\end{equation}
By the fact that $S^{1}$ is a $K(\Gamma/F_{2},1)$ and the above calculations, we see the spectral sequence above degenerates at the second page because $H^{p}(\Gamma/F_{2};H^{q}(F_{2};\slf(4,\R))) = 0$ for all $p > 1$.  Thus we have the isomorphisms in Equation \ref{eqisos}.
\begin{align}
H^{0}(\Gamma/F_{2};H^{1}(F_{2};\slf(4,\R))) &\simeq H^{1}(\Gamma; \slf(4,\R))  \label{eqisos}\\
 H^{1}(\Gamma/F_{2};H^{1}(F_{2};\slf(4,\R))) &\simeq H^{2}(\Gamma; \slf(4,\R)) \nonumber
\end{align}
Note that $H^{0}(\Gamma; \slf(4,\R)) = 0$ for similar reasons established above, whereas $H^{3}(\Gamma; \slf(4,\R)) = 0$ by a Poincar\'{e} duality argument \cite{Heusener2011Infinitesimal}.  Thus Equations \ref{eqisos} provides us a means of expressing the cohomology of $\Gamma$ entirely in terms of the invariants and co-invariants of the cyclic $\Gamma/F_{2}$-action on $H^{1}(F_{2};\slf(4,\R))$.  Consequently if we are able to understand the dynamics of the $\Gamma/F_{2}$-action on $H^{1}(F_{2};\slf(4,\R))$, we are able to understand the cohomology of the $\Gamma$-module $\slf(4,\R)$.
\end{subsection}

\begin{subsection}{Induced Action on Tangent Space of the Character Variety}\label{ssinducedaction}
Before proceeding we rephrase the $\Gamma/F_{2}$-action in terms of the differential of an action on the entire character variety of $F_{2}$ into $\PGL(4,\R)$.  We denote $\Hom(F_{2},\PGL(4,\R))$ by the representation variety of $F_{2}$ into $\PGL(4,\R)$.  The representation variety of $F_{2}$ enjoys a rich action via the entire outer automorphism group of $F_{2}$, $\text{Out}(F_{2})$.  Specifically, the automorphism induced by $\alpha \in \text{Out}(F_{2})$ takes a representation $\sigma: F_{2} \lra \PGL(4,\R)$ and sends it to $(\alpha\sigma): F_{2} \lra \PGL(4,\R)$ defined by $(\alpha\sigma)(g) := \sigma(\alpha^{-1}g)$. \\
\\
In our context, we have $\Gamma$ is a hyperbolic once-punctured torus bundle with monodromy $\phi$ and holonomy representation $\rho: \Gamma \lra \PSO(3,1) \subset \PGL(4,\R)$.  The monodromy $\phi$ corresponds to a homeomorphism of the once-punctured torus fiber and thus induces an action on its fundamental group which we abusively denote by $\phi$ as well.  This automorphism may be identified with an element of the outer automorphism group of $F_{2}$ as explained in Section \ref{ssoptb}, thus we have an automorphism of the representation variety induced by $\phi \in \Out(F_{2})$.  \\
\\
Recall that $\Gamma$ admits a presentation of the form in Equation \ref{eqpres} where conjugation by $x$ on the free group generators $a,b$ acts by the outer automorphism $\phi$.  This means that for any element $g \in \Gamma$, we may post-compose the action of $\phi \in \Out(F_{2})$ with conjugation by $\rho(g)$.  Applying this action with $g = x$ to the generators of an arbitrary $\sigma \in \Hom(F_{2},\PGL(4,\R))$ yields the equations below.
\begin{align*}
\rho(x)(\phi\sigma)(a)\rho(x)^{-1} &= \rho(x)\sigma\left(\phi^{-1}(a)\right)\rho(x)^{-1} \\
\rho(x)(\phi\sigma)(b)\rho(x)^{-1} &= \rho(x)\sigma\left(\phi^{-1}(b)\right)\rho(x)^{-1} 
\end{align*}
Note that when $\sigma = \rho$, we have the above equations reduce to the identities below. 
\begin{align*}
\rho(x)(\phi\rho)(a)\rho(x)^{-1} &= \rho(a) \\
 \rho(x)(\phi\rho)(b)\rho(x)^{-1} &= \rho(b) 
\end{align*}
That is to say, this action has the holonomy representation of $\rho \in \Hom(F_{2},\PGL(4,\R))$ as a \emph{fixed point}.  In addition, if we choose different choices of $g$ inside the coset $xF_{2}$, we result with conjugate actions.  Thus, this action descends to an action on the \emph{character variety} of $F_{2}$ into $\PGL(4,\R)$, denoted by $\chi(F_{2},\PGL(4,\R))$.  Because the action restricted to $F_{2}$ acts \emph{trivially} on the character variety, we have an induced $\Gamma/F_{2}$-action on the character variety.  We abusively denote this representation by $\phi$. \\ 
\\
This action relates to the action of $\Gamma/F_{2}$ on $H^{1}(F_{2};\slf(4,\R))$ in the following manner.  The holonomy representation of $\rho: F_{2} \lra \PGL(4,\R)$ is a fixed point in the character variety of $F_{2}$ into $\PGL(4,\R)$ under the action $\phi$.  One may think of the $\Gamma/F_{2}$-action on $H^{1}(F_{2};\slf(4,\R))$ as the \emph{differential} of the action of $\phi$ at the point $[\rho] \in \chi(F_{2},\PGL(4,\R))$.  First we explain why the point $[\rho] \in \chi(F_{2},\PGL(4,\R))$ should be a smooth point of the character variety.  In the representation variety, $\Hom(F_{2}, \PGL(4,\R))$ admits a particularly nice smooth structure as $\PGL(4,\R)\times \PGL(4,\R)$ once we choose a generating set for $F_{2}$.  Thus we have that $\rho$ corresponds to two matrices $\rho(a)$ and $\rho(b)$ which represent the fundamental group in $\PSO(3,1)$ of the once-punctured fiber.  Because $\rho$ is irreducible, so too are all nearby representations, thus they have trivial stabilizers in $\PGL(4,\R)$ \cite{Ballas2018Convex}.  Consequently $[\rho]$ is a smooth point in the character variety of $\chi(F_{2},\PGL(4,\R))$.  Because $[\rho]$ is a smooth point we may identify the tangent space of $\chi(F_{2},\PGL(4,\R))$ at the geometric representation $[\rho]$ with $H^{1}(F_{2};\slf(4,\R))$ \cite{sikora}.  \\
\\
Let $\rho_{t}$ be a smooth path of representations in $\Hom(F_{2},\PGL(4,\R))$ so that $\rho_{0} = \rho$ and let $\mathfrak{g}$ denote the $F_{2}$-module $\slf(4,\R)_{\Ad \rho}$.  In the typical fashion, this determines a crossed homomorphism $f: F_{2} \lra \g$ as described in Section \ref{subsec-twistedcohom}.
\begin{equation*}\label{eqcocycle}
f(g) := \frac{d}{dt}\bigg|_{t=0} \rho_{t}(g)\rho^{-1}(g) 
\end{equation*}
The crossed homomorphism satisfies the property that $f(gh) = f(g) + gf(h)$ for all $h,g \in F_{2}$ where here left multiplication by $g$ denotes the $\Gamma$-action on $\mathfrak{g}$ restricted to $F_{2}$.  This correspondence determines an isomorphism between the tangent space of the representation variety of $F_{2}$ at the point $\rho$ and the $1$-cocycles of the group cohomology determined by the $\Gamma$-module $\mathfrak{g}$.  Note the $1$-cocycles of $F_{2}$ are entirely determined by two choices of vectors $X, Y \in \mathfrak{g}$, thus $Z^{1}(F_{2}; \mathfrak{g}) \simeq \mathfrak{g}\times \mathfrak{g}$.  This can also be seen simply by identifying $\rho_{t}$ with a smooth path of matrices $a_{t}$ and $b_{t}$ so that $X$ and $Y$ are the derivatives of $a_{t}a_{0}^{-1}$ and $b_{t}b_{0}^{-1}$ respectively.  \\
\\
We now study how the action of $\phi$ on the representation variety affects a cocycle in $Z^{1}(F_{2}; \mathfrak{g})$.  Recall $\phi$ acts on a representation $\sigma$ by first pre-composing it with the inverse of the outer-automorphism $\phi$, then conjugating by $\rho(x)$.  Notationally, our smooth path $\rho_{t}$ will get sent to $(\phi\rho_{t})$.  We denote the associated cocycle by $(\phi f)$ which we evaluate in Equation \ref{eqactedcocycle} below.  Taking the derivative of $(\phi\rho_{t})$ evaluated at $g \in F_{2}$ yields
\begin{align}
(\phi f)(g) &:= \frac{d}{dt}\bigg|_{t=0} \left[(\phi\rho_{t})(g)\right]\left[(\phi\rho)(g)\right]^{-1} \nonumber \\
&= \frac{d}{dt}\bigg|_{t=0} \left[\rho(x)\rho_{t}\left(\phi^{-1}(g)\right)\rho(x)^{-1} \right]  \left[\rho(x)\rho\left(\phi^{-1}(g)\right)\rho(x)^{-1} \right]^{-1} \nonumber \\
&= \frac{d}{dt}\bigg|_{t=0} \left[\rho(x)\rho_{t}\left(\phi^{-1}(g)\right) \rho\left(\phi^{-1}(g)\right)^{-1}\rho(x)^{-1} \right]  \nonumber \\
&=\rho(x) \frac{d}{dt}\bigg|_{t=0} \rho_{t}\left(\phi^{-1}(g)\right) \rho\left(\phi^{-1}(g)\right)^{-1} \rho(x)^{-1} \nonumber \\
&= \rho(x)f\left(\phi^{-1}(g)\right)\rho(x)^{-1} \nonumber\\
&= \rho(x)f\left(x^{-1}gx\right)\rho(x)^{-1} \label{eqactedcocycle}
\end{align}
Adopting the notation of the $\Gamma$-action on $\mathfrak{g}$, this is to say that $(\phi f)(g) = xf(x^{-1}gx)$ for all $g \in F_{2}$.  One can readily verify in the case that $f$ is a coboundary as defined in Section \ref{ssiaata}, that $(\phi f)(g) = (1-g)(xY)$ for some $Y \in \mathfrak{g}$.  Thus $\phi$ sends coboundaries to coboundaries, and the action on $Z^{1}(F_{2};\mathfrak{g})$ descends to an action on $H^{1}(F_{2};\mathfrak{g})$.  In Equation \ref{eqactedcocycle} if we instead of conjugate by $\rho(h)$ where $h \in F_{2}$ instead of $\rho(x)$, we may expand $f(h^{-1}gh)$ via the properties of crossed homomorphisms as $h^{-1}gh \in F_{2}$.  This yields the equation below.
\begin{equation*}
\rho(h)f(h^{-1}gh)\rho(h)^{-1} = h\left(\left(h^{-1}g-h^{-1}\right)f(h) + h^{-1}f(g)\right)= f(g) + (g-1)f(h) 
\end{equation*}
Thus we have that the action by $F_{2}$ acts trivially on $H^{1}(F_{2};\mathfrak{g})$ and we recover the induced action of $\Gamma/F_{2}$ on $H^{1}(F_{2};\mathfrak{g})$ arising from the Lyndon-Hochschild-Serre spectral sequence that is expressed entirely from a group theoretic perspective.  We summarize these results in the finding below.
\begin{theorem}\label{thm1p5}
Let $M$ be a hyperbolic once-punctured torus bundle with holonomy $\rho: \Gamma \lra \PSO(3,1) \subset \PGL(4,\R)$ and monodromy $\phi \in \emph{\Out}(F_{2})$.  The equivalence class of the holonomy representation restricted to the fundamental group of a once-punctured torus fiber is a smooth point of the character variety of $F_{2}$ into $\PGL(4,\R)$.  Moreover, $[\rho]$ is a fixed point in the character variety under the $\Gamma/F_{2}$-action induced by $\phi$.  The differential of the $\Gamma/F_{2}$-action at this fixed point is equal to the group theoretic action of the Lyndon–Hochschild–Serre spectral sequence associated to the $\Gamma$-module $\slf(4,\R)_{\Ad \rho}$ and its normal subgroup $F_{2}$.  
\end{theorem}

This result relates the dynamical action of the outer automorphism of $\phi$ acting on the character variety and the purely group theoretic action arising from cohomology.  By Equation \ref{eqisos}, to determine whether or not a hyperbolic once-punctured  torus bundle is infinitesimally projectively rigid rel cusp, it suffices to understand the dynamics of this action at the geometric representation in the character variety, specifically, how the differential acts on the tangent space.  
\end{subsection}

\begin{subsection}{Induced Action and Twisted Alexander Polynomial}\label{ssiaata}
In this section we study one of the most natural invariants of the $\Gamma/F_{2}$-action on $H^{1}(\Gamma; \mathfrak{g})$: its characteristic polynomial.  We show that it is in fact equal to the twisted Alexander polynomial of the representation $\Ad \rho: \Gamma \lra \Aut(\mathfrak{g})$.  To this end let $\Gamma$ be the fundamental group of a hyperbolic once-punctured torus bundle and have the presentation as in Equation \ref{eqpres}.  Pick $a,b$ as a choice of generating set for $F_{2}$ and the relations imposed by $\Gamma$ are given by Equation \ref{eqrelations} below.  
\begin{equation}\label{eqrelations}
R_{1} = \phi(a)xa^{-1}x^{-1},\phantom{==}R_{2} = \phi(b)xb^{-1}x^{-1}
\end{equation}
Let $\alpha: \Gamma \lra \langle t \rangle$ be the map that takes $F_{2}$ to $1$ and $x$ to $t$.  Following the calculations introduced in Section \ref{subsec-twistedalex} we see that $\Phi: \Z F_{3} \lra \Aut(\mathfrak{g})$ takes $a$ to $\Ad\rho(a)$, $b$ to $\Ad\rho(b)$, and $x$ to $t\Ad\rho(x)$.  Instead of taking the taking Fox-derivatives in the variables $a$, $b$ and $x$ above in Equation \ref{eqrelations}, we take the derivative of equivalent relations in terms of $\phi^{-1}$, as $\phi^{-1}$ will appear in the characteristic polynomial of the $\Gamma/F_{2}$-action.  An equivalent way of writing $R_{1}$ and $R_{2}$ is
\begin{equation*}
R_{1}' = \phi^{-1}(a)x^{-1}a^{-1}x,\phantom{==}R_{2}' = \phi^{-1}(b)x^{-1}b^{-1}x
\end{equation*}
We emphasize that $\phi^{-1}(g)$ denotes the \emph{functional inverse} of $\phi$ applied to $g$, not the element $\phi(g)^{-1}$.  Taking the Fox-derivatives of $R_{1}'$ yields the equations below.  
\begin{align*}
\frac{\del}{\del a}R_{1}' &= \frac{\del}{\del a} \left(\phi^{-1}(a)\right) + \phi^{-1}(a) \frac{\del}{\del a} \left(x^{-1}a^{-1}x\right) \nonumber \\
&=  \frac{\del}{\del a} \left(\phi^{-1}(a)\right)  - \phi^{-1}(a)x^{-1}a^{-1} = \frac{\del}{\del a} \left(\phi^{-1}(a)\right)  - x^{-1} \\
\frac{\del}{\del b}R_{1}' &= \frac{\del}{\del b} \left(\phi^{-1}(a)\right) + \phi^{-1}(a) \frac{\del}{\del b} \left(x^{-1}a^{-1}x\right) \nonumber \\
&=  \frac{\del}{\del b} \left(\phi^{-1}(a)\right) \\
\frac{\del}{\del x}R_{1}' &= \frac{\del}{\del x} \left(\phi^{-1}(a)\right) + \phi^{-1}(a) \frac{\del}{\del x} \left(x^{-1}a^{-1}x\right) \nonumber \\
&= \phi^{-1}(a)(x^{-1}a^{-1}-x^{-1}) = x^{-1}-x^{-1}a
\end{align*}
An entirely analogous calculation for $R_{2}'$ may be carried.  Applying $\Phi$ yields the twisted Alexander matrix and taking the minor associated to the 3rd column yields the twisted Alexander polynomial below in Equation \ref{eqtwistedminorcalc}.  Note the denominator is non-zero as $x$ is sent to the non-trivial element $t$ under $\alpha$.  
\begin{footnotesize}
\begin{equation}\label{eqtwistedminorcalc}
\Delta_{\rho,\alpha}(t) = \det \left( \begin{array}{ccc}
\Ad \rho\left(\frac{\del}{\del a} \left(\phi^{-1}(a)\right) \right) - t\Ad \rho(x^{-1})  & \Ad \rho\left(\frac{\del}{\del b} \left(\phi^{-1}(a)\right)\right)   \\
\Ad \rho\left(\frac{\del}{\del a} \left(\phi^{-1}(b)\right)\right)    & \Ad \rho\left(\frac{\del}{\del b} \left(\phi^{-1}(b)\right)\right) - t\Ad \rho^{-1}(x) 
\end{array}\right)/\det(\text{Id}_{\g}-t\Ad \rho(x))
\end{equation}
\end{footnotesize}
Let us now instead compare this to the characteristic polynomial of the $\Gamma/F_{2}$-action on $H^{1}(F_{2};\g)$.  As a preliminary observation note that the $\phi$-action on $\g \times \g$ as defined in Equation \ref{eqactedcocycle} descends to a generator of the $\Gamma/F_{2}$-action on $H^{1}(F_{2};\mathfrak{g})$.  Using Equation \ref{eqactedcocycle}, the $\phi$-action takes a cocycle $f : F_{2} \lra \g$ and sends it to
\begin{equation*}
(\phi f)(g) = \rho(x)f(x^{-1}gx)\rho(x)^{-1} = xf(x^{-1}gx) = xf(\phi^{-1}g)
\end{equation*}
where the last equality is written in terms of the $\Gamma$-action on $\mathfrak{g}$.  Because every cocycle is determined by its value on $a$ and $b$, $Z^{1}(F_{2};\mathfrak{g})$ is isomorphic to $\g\times \g$ as vector spaces.  Thus $\phi f$ is determined by its action on $a$ and $b$.  Writing this out explicitly using the properties of crossed homomorphisms yields
\begin{align}
(\phi f)(a) &= xf(x^{-1}ax) = xf\left(\phi^{-1}(a)\right) = x\left(\frac{\del}{\del a}\phi^{-1}(a)\right)f(a) + x\left(\frac{\del}{\del b}\phi^{-1}(a)\right)f(b) \nonumber \\
(\phi f)(b) &= xf(x^{-1}bx) = xf\left(\phi^{-1}(b)\right) = x\left(\frac{\del}{\del b}\phi^{-1}(b)\right)f(a) + x\left(\frac{\del}{\del b}\phi^{-1}(b)\right)f(b) \label{eqcocycleaction2}
\end{align}
This means the characteristic polynomial $p_{\g \times\g}(t)$ of $\phi$ acting on $\mathfrak{g}\times \mathfrak{g}$ in the variable $t$ is given by
\begin{equation}\label{eqEQACTIONgxg2}
p_{\g\times\g}(t) = \det \left(
\begin{array}{cc}
\Ad\rho \left(x\frac{\del}{\del a}\left(\phi^{-1}(a)\right)\right) - t\cdot\text{id}_{\g} &  \Ad \rho \left(x\frac{\del}{\del b}\left(\phi^{-1}(a)\right)\right)\\
\Ad \rho\left( x\frac{\del}{\del b}\left(\phi^{-1}(b)\right) \right) &  \Ad \rho \left(x\frac{\del}{\del b}\left(\phi^{-1}(b)\right) \right)- t\cdot\text{id}_{\g} 
\end{array}
\right)
\end{equation} 
Notice in Equation \ref{eqEQACTIONgxg2}, there is a common term of $\Ad \rho(x)$.  If we factor this out Equation \ref{eqEQACTIONgxg2} becomes 
\begin{equation*}\label{eqcharpolyreduced}
\det \Ad\rho(x) \det \left(
\begin{array}{cc}
\Ad\rho \left(\frac{\del}{\del a}\left(\phi^{-1}(a)\right)\right) - t\Ad\rho(x^{-1}) &  \Ad \rho \left(\frac{\del}{\del b}\left(\phi^{-1}(a)\right)\right)\\
\Ad \rho\left( \frac{\del}{\del b}\left(\phi^{-1}(b)\right) \right) &  \Ad \rho \left(\frac{\del}{\del b}\left(\phi^{-1}(b)\right) \right)- t\Ad\rho(x^{-1}) 
\end{array}
\right)
\end{equation*} 
Because the adjoint representation has determinant one, the equation above is readily seen to be the numerator of Equation \ref{eqtwistedminorcalc}.  To calculate the characteristic polynomial of the $\Gamma/F_{2}$-action on $H^{1}(F_{2};\g)$ we recall the short exact sequence of $\Gamma$-modules given below in Equation \ref{eqsesgammods}.  
\begin{equation}\label{eqsesgammods}
B^{1}(F_{2};\mathfrak{g}) \lra Z^{1}(F_{2};\mathfrak{g}) \lra H^{1}(F_{2};\mathfrak{g})
\end{equation}
We know from our previous analysis $Z^{1}(F_{2};\mathfrak{g}) \simeq \g \times \g$ where the $\phi$-action on $Z^{1}(F_{2};\mathfrak{g})$ is defined Equation \ref{eqactedcocycle}.  We now consider how the $\phi$-action acts on the subspace of coboundaries $B^{1}(F_{2};\g) \subset Z^{1}(F_{2}; \g)$.  As a vector space $B^{1}(F_{2};\g)$ is all crossed homomorphisms $f: \Gamma \lra \g$ satisfying $f(g) = (1-g)X$ for some $X \in \g$ and for all $g \in G$.  Because the $F_{2}$-action on $\g$ has no non-trivial invariants, the natural map taking $\g \lra B^{1}(F_{2};\g)$ via $X$ gets sent to the coboundary $f(g) = (1-g)X$ is an isomorphism of vector spaces.    \\
\\
With this in mind we calculate the $\phi$-action on $B^{1}(F_{2};\g)$ via the isomorphism $\g \simeq B^{1}(F_{2};\g)$.  
\begin{align*}
(\phi f)(g) &= \rho(x)f\left(x^{-1}gx\right)\rho(x)^{-1} = \rho(x)\left(1-x^{-1}gx\right)X\rho(x)^{-1} = x(1-x^{-1}gx)X \\
&= (1 - g)(xX) 
\end{align*}
Thus, the $\phi$-action on $\g$ is simply the adjoint representation taking $X \in \g$ to $xX = \rho(x)X\rho(x)^{-1}$.  Therefore the characteristic polynomial of the $\Gamma/F_{2}$-action is the quotient of the two $\phi$-actions on the cocycles and coboundaries in Equation \ref{eqsesgammods}.  The characteristic polynomial of $\phi$ on the cocycles, $p_{\g\times \g}$, is given by Equation \ref{eqEQACTIONgxg2} whereas the characteristic polynomial of $\phi$ on the coboundaries, $p_{\g}$ is given by $\det(\Ad\rho(x) - t\text{Id}_{\g})$.  Using the fact that $\Ad \rho$ preserves the Killing Form, we know that $\Ad \rho(x)$ and $\Ad \rho(x^{-1})$ are similar and thus share the same characteristic polynomial.  This means that
\begin{equation*}
\det(\Ad\rho(x) - t\text{Id}_{\g}) = \det\Ad\rho (x) \det \left(\text{Id}_{\g} - t \Ad \rho(x)^{-1}\right) = \det \left(\text{Id}_{\g} - t \Ad \rho(x)\right) 
\end{equation*}
Thus we have the following theorem.
\begin{theorem}\label{thm1}
Let $\emph{\Ad} \rho : \Gamma \lra \emph{\Aut}(\g)$ be the composition of the holonomy representation of $\Gamma \lra \emph{\PSO}(3,1) \subset \emph{\PGL}(4,\R)$ of a once-punctured hyperbolic torus bundle $M$ and the adjoint representation of $\emph{\PGL}(4,\R)$.  Let the monodromy of $M$ be given by $\phi \in \emph{\SL}(2,\Z)$.  Express $\Gamma = \langle a,b,x \, |\,  \phi(a) = xax^{-1} \text{ and } \phi(b) = xbx^{-1}\rangle$.  Then the characteristic polynomial the $\Gamma/F_{2}$-action induced by the Lyndon–Hochschild–Serre spectral sequence on group cohomology is equal to the twisted Alexander polynomial of the representation $\emph{\Ad} \rho$ with the map $\alpha : \Gamma \lra \Gamma/F_{2} \simeq \langle x \rangle$.  
\end{theorem}
In the next section we will develop tools via cohomology calculations to detect whether or not a hyperbolic once-punctured torus bundle is infinitesimally projectively rigid rel cusp in terms of the roots of the twisted Alexander polynomial which by Theorem \ref{thm1} are precisely the eigenvalues of the $\Gamma/F_{2}$-action on the tangent space, $H^{1}(F_{2};\g)$, to the character variety at the geometric representation $[\rho]$.  
\end{subsection}

\begin{subsection}{Cohomology of the Once-Punctured Torus}\label{sscohomopt}
We begin by observing we can calculate the dimension of the cohomology of the once-punctured torus at the geometric representation $\rho: \Gamma \lra \PSO(3,1) \subset \PGL(4,\R)$ by simple dimension counting.  Referring to Equation \ref{eqsesgammods}, we see that 
\begin{equation*}
\dim H^{1}(F_{2};\g) = \dim Z^{1}(F_{2};\g) - \dim B^{1}(F_{2};\g) = \dim G^{2} - \dim G = \dim G
\end{equation*}
The quantity we are mostly interested in are the $\Gamma/F_{2}$-invariants of $H^{1}(F_{2};\g)$.  As was done in Section \ref{sscohomopt}, we will let $x$ be the meridian whose image under the quotient map $\Gamma \lra \Gamma/F_{2}$ generates the quotient group.  The meridian $x$ commutes with the longitude $l = aba^{-1}a^{-1}$.  Their image under the holonomy map forms a generating set for the fundamental group of the cusp of our hyperbolic once-punctured torus bundle.  \\
\\
Consider the good pair $(\del F, F)$ where $\del F$ is the longitude of the once-punctured torus fiber $F$.  We may form the long-exact sequence in cohomology associated to this pair below in Equation \ref{eqlongexactgp}.
\begin{small}
\begin{equation}\label{eqlongexactgp}
H^{0}(F, \del F) \rightarrow H^{0}(F) \rightarrow H^{0}(\del F) \rightarrow H^{1}(F, \del F) \rightarrow H^{1}(F) \rightarrow H^{1}(\del F) \rightarrow H^{2}(F, \del F) \rightarrow 0
\end{equation}
\end{small}
We are suppressing the coefficient module $\mathfrak{g}_{\Ad \rho}$ in Equation \ref{eqlongexactgp} to save space.  Note for $k = 0,1,2$, we have the perfect pairing 
\begin{equation}\label{eqperfpair}
H^{k}(F,\del F, \g)\otimes H^{2-k}(F, \g) \lra H^{2}(F,\del F, \g\otimes \g) \lra H^{2}(F,\del F, \R)
\end{equation}
where the last map on the coefficient module is given by the non-degenerate $\Gamma$-invariant bilinear Killing form $B: \g\otimes\g \lra \R$.  This pairing is $\Gamma/F_{2}$-invariant, and thus as $\Gamma/F_{2}$-modules we have that $H^{1}(F)$ and $H^{1}(F,\del F)^{*}$ are isomorphic as $\Gamma/F_{2}$-modules.  Since $\dim H^{0}(F) = 0$, as there are no $F_{2}$-invariants in $\g$, this means $\dim H^{2}(F,\del F) = 0$.  Thus the exact sequence in Equation \ref{eqlongexactgp} reduces to Equation \ref{eqgoodpair} below.  
\begin{equation}\label{eqgoodpair}
0 \lra H^{0}(\del F) \lra H^{1}(F, \del F) \lra H^{1}(F) \lra H^{1}(\del F) \lra 0
\end{equation}
The image of the longitude $\del F$ is parabolic under $\rho$, so up to conjugation $\rho(l) = \rho(aba^{-1}b^{-1})$ is a unipotent matrix whose Jordan decomposition is of the form below \cite{Ballas2014Deformations}.
\begin{equation}\label{eqlong}
\rho(l) = \left(
\begin{array}{cccc}
1 & 0 & 0 & 0 \\
0 & 1 & 1 & 0 \\
0 & 0 & 1 & 1 \\
0 & 0 & 0 & 1
\end{array}
\right)
\end{equation}
A straight forward calculation determines the infinitesimal centralizer of $\rho(l)$ in $\slf(4,\R)$, namely the $X \in \slf(4,\R)$ for which $\rho(l)X\rho(l)^{-1} = X$, is 5-dimensional.  A 2-dimensional subspace of the infinitesimal centralizer lies inside $\pso(3,1)$ and the remaining 3-dimensions lie inside $\vv$ as defined in Section \ref{sscohomology}.  \\
\\
By Equation \ref{eqperfpair}, we know that $H^{1}(F)$ and $H^{1}(F,\del F)^{*}$ are isomorphic as $\Gamma/F_{2}$-modules.  Thus the characteristic polynomials of their $\Gamma/F_{2}$-actions are the same.  Moreover, the characteristic polynomial of the $\Gamma/F_{2}$-action on $H^{1}(F,\del F)^{*}$ is equal to the characteristic polynomial of the $\Gamma/F_{2}$-action on $H^{1}(F,\del F)$.  Thus, if we wish to calculate the characteristic polynomial of the $\Gamma/F_{2}$-action on $H^{1}(F)$ it suffices instead to calculate it on $H^{1}(F,\del F)$.  This is for computational purposes easier insofar as the corresponding chain complex associated to $\del F \subset F$ is given by the commutative diagram below.
\begin{equation}\label{eqchaincohom}
\begin{tikzcd}
\ker \res_{l}^{B} \arrow{r}{i} \arrow{d}{i}			&B^{1}(F_{2};\g) \arrow{r}{\res_{l}^{B}}   \arrow{d}{i}				& B^{1}(\Z;\g) \arrow{d}{i}\\
\ker \res_{l}^{Z} \arrow{r}{i}\arrow{d}{q} 	&Z^{1}(F_{2};\g) \arrow{r}{\res_{l}^{Z}}	\arrow{d}{q}			& Z^{1}(\Z;\g) \arrow{d}{q}\\
 H^{1}(F,\del F; \g) \arrow{r}{i}		 	& H^{1}(F_{2};\g) \arrow{r}{\res_{l}}							& H^{1}(\Z;\g) \\
\end{tikzcd}
\end{equation}
In the above diagram, the $\Z$ on the right-hand column is generated by $l = aba^{-1}b^{-1}$.  Because the restriction $\res_{l}^{B}$ is an isomorphism, which up to identification is simply the identity map on $\g$, the top-left vertex is trivial.  Consequently the relative cohomology group $H^{1}(F,\del F; \g)$ is isomorphic to the kernel of $\res_{l}^{Z}: Z^{1}(F_{2};\g) \lra Z^{1}(\Z;\g)$. \\
\\ 
Up to isomorphism, we know that $Z^{1}(F_{2};\g)$ is $\g\times \g$, thus we may think of $H^{1}(F,\del F; \g)$ as the subspace of $\g\times\g$ for which the deformations in these directions produce no change in the holonomy of $l$.  That is to say, these are the directions that deform while still preserving the longitude as a parabolic element of the form in Equation \ref{eqlong}.  Using Fox-derivatives to calculate $\res_{l}^{Z}(X,Y)$ we see that
\begin{equation}\label{eqresl}
\res_{l}^{Z}(X,Y) = (1-aba^{-1})X + (a-aba^{-1}b^{-1})Y
\end{equation}
So if we have the explicit representation $\Ad \rho: \Gamma \lra \Aut(\g)$, we can explicitly calculate $H^{1}(F,\del F; \g)$, and correspondingly, the characteristic polynomial of the $\Gamma/F_{2}$-action.  \\
\\
The first part of Equation \ref{eqgoodpair} yields the short exact sequence
\begin{equation}\label{eqsesker}
0 \lra H^{0}(\del F) \lra H^{1}(F, \del F) \lra K \lra 0
\end{equation}
where $K$ is the kernel of the map $\res_{l} : H^{1}(F) \lra H^{1}(\del F)$.  Thus the characteristic polynomial of $\Gamma/F_{2}$ on $H^{1}(F,\del F)$ splits as the characteristic polynomial on the pieces $H^{0}(\del F)$ and $K$.  The $\Gamma/F_{2}$-action on $H^{0}(\del F)$ is up to isomorphism simply conjugation by $\rho(x)$ on the invariants of $\rho(l)$ and thus $\Gamma/F_{2}$ acts unipotently on $H^{0}(\del F)$ as $\rho(x)$ is unipotent.  Because this is a 5-dimensional space, this means the characteristic polynomial of $H^{1}(F,\del F)$ has a factor of $(1-t)^{5}$. \\
\\
However if this is the largest factor, this means there are no $\Gamma/F_{2}$-invariants in $K$.  If we take the corresponding long exact sequence of $\Gamma/F_{2}$-cohomology of the short-exact sequence of $\Gamma/F_{2}$-modules in Equation \ref{eqsesker}, we obtain the exact sequence below.
\begin{small}
\begin{equation}\label{eqcohomshortlongexact}
H^{0}(\del F)^{\Gamma/F_{2}} \lra H^{1}(F,\del F)^{\Gamma/F_{2}} \lra K^{\Gamma/F_{2}} \lra H^{0}(\del F)_{\Gamma/F_{2}} \lra H^{1}(F,\del F)_{\Gamma/F_{2}} \lra K_{\Gamma/F_{2}} 
\end{equation}
\end{small}
Because there are no $\Gamma/F_{2}$-invariants in $K$, this means that $H^{0}(\del F)^{\Gamma/F_{2}} \lra H^{1}(F,\del F)^{\Gamma/F_{2}}$ in as isomorphism by exactness of the above sequence.  Moreover, because the $\Gamma/F_{2}$-action on $H^{0}(\del F)$ is isomorphic to the $x$-action on $H^{0}(\del F)$ we see that 
\begin{equation*}
H^{0}(\del F)^{\Gamma/F_{2}} \simeq H^{0}(\del F)^{\Ad \rho \, \pi_{1}(\del M)}
\end{equation*}
A straight-forward calculation of the infinitesimal stabilizer shows that this means $\dim H^{0}(\del F)^{\Gamma/F_{2}}$ $= \dim H^{0}(\del F)^{\Ad \rho \, \pi_{1}(\del M)} = 3$ \cite{Ballas2018Convex}.  Thus $M$ is infinitesimally projectively rigid.  We summarize these findings in the following theorem.
\begin{theorem}\label{thm2}
Let $M$ be a hyperbolic once-punctured torus bundle with holonomy $\rho: \Gamma \lra \emph{\PSO}(3,1) \subset \emph{\PGL}(4,\R)$.  Let $\g$ denote the $\Gamma$-module $\slf(4,\R)$ associated to the composition $\emph{\Ad} \rho: \Gamma \lra \emph{\Aut}(\g)$.  If the twisted Alexander polynomial of $\emph{\Ad} \rho$ has 1 as root of multiplicity 5, then $M$ is infinitesimally projectively rigid rel cusp.  
\end{theorem}
For computational purposes, we can speed up calculations by the following observation.  By Equation \ref{eqhsplit}, $H^{1}(\Gamma;\g)^{\Gamma/F_{2}}$ splits into the direct sum of $H^{1}(\Gamma;\pso(3,1))^{\Gamma/F_{2}}$ and  $H^{1}(\Gamma;\vv)^{\Gamma/F_{2}}$.  By Thurston's Dehn-Filling theorem, we know that $\dim H^{1}(\Gamma;\pso(3,1))^{\Gamma/F_{2}} = 2$, thus to verify infinitesimal projective rigidity we may carry out the calculations of Equation \ref{eqsesker} and Equation \ref{eqcohomshortlongexact} with coefficients in the subspace $\vv$ instead of the full $\slf(4,\R)$.  As remarked earlier, in this case $\dim H^{0}(\del F) = 3$, thus, we have the following corollary. 
\begin{corollary}
Let $M$ be a hyperbolic once-punctured torus bundle with holonomy $\rho: \Gamma \lra \emph{\PSO}(3,1) \subset \emph{\PGL}(4,\R)$.  Let $\g$ denote the $\Gamma$-module associated to the composition $\emph{\Ad} \rho: \Gamma \lra \emph{\Aut}(\g)$ and $\vv$ denote the orthogonal complement of $\pso(3,1)$ under the Killing form of $\g$.  If the twisted Alexander polynomial of $\emph{\Ad} \rho$ has 1 as root of multiplicity 3, then $M$ is infinitesimally projectively rigid rel cusp. 
\end{corollary}
If one wishes to proceed by calculating the twisted Alexander polynomial, this is a bit easier in calculations, as $\vv$ is 9-dimensional whereas $\slf(4,\R)$ is 15-dimensional.  Thus, $H^{1}(F,\del F; \vv)$ is 9-dimensional, and the characteristic polynomial of the $\Gamma/F_{2}$-action is degree 9 instead of degree 15, so it will be easier to factor or implement some root isolation algorithm.  
\end{subsection}

\begin{subsection}{Yet Another Twisted Polynomial Calculation}
We conclude this section with a way to calculate the twisted Alexander polynomial closer in the spirit to its definition than the characteristic polynomial of the $\Gamma/F_{2}$-action on $H^{1}(F_{2};\g)$.  To this end let $M$ be a hyperbolic once-punctured torus bundle with fundamental group $\Gamma$ and holonomy $\rho: \Gamma \lra \PSO(3,1) \subset \PGL(4,\R)$.  Let us lift the representation into $\SO(3,1)$ which can be done by Culler \cite{CULLER198664}, and denote the $\Gamma$-module $\R^{3,1}$ by $V$ where $\Gamma$ acts on $V$ via the so-called `birth certificate representation' $gv = \rho(g)v$ for each $v \in V$ and $g \in G$.  We are able to relate this to the adjoint representation in the following manner.  Due to the $\Gamma$-invariant biliear form as defined in Equation \ref{eq-lorentz}, we have that $V$ is isomorphic to $V^{*}$ as $\Gamma$-modules.  \\
\\
As $\Gamma$-modules, $V^{*}\otimes V$ is naturally isomorphic to $\Hom(V,V)$ via the isomorphism taking $f\otimes v$ to the linear map $L(w):=f(w)v$.  If we carefully follow the $\Gamma$-action on $V^{*}\otimes V$ through this isomorphism, we find that the induced $\Gamma$-action on $\Hom(V,V)$ is simply conjugation.  Thus, if we pick a basis for $V$ we may identify $\Hom(V,V)$ with $\mathfrak{gl}(4,\R)$.  The image of our representation $\Gamma \lra \Aut(V^{*}\otimes V)$ preserves the invariant subspace of \emph{traceless} matrices, i.e. $\slf(4,\R) \subset \mathfrak{gl}(4,\R)$.  Thus we have the short-exact sequence of $\Gamma$-modules below.
\begin{equation*}
\mathfrak{sl}(4,\R) \lra \mathfrak{gl}(4,\R) \lra \R
\end{equation*}
where the quotient inherits a trivial $\Gamma$-action as $\det \circ \rho : \Gamma \lra \R$ is trivial.  This short exact sequence yields that twisted Alexander polynomial of the representation of $\tilde{\rho}\otimes\tilde{\rho}$ into $\mathfrak{gl}(4,\R)$ will be the product of the twisted Alexander polynomials of $\mathfrak{sl}(4,\R)$ and $\R$.  Symbolically that is to say
\begin{equation*}
\Delta_{\mathfrak{gl}}(t) = \Delta_{\mathfrak{sl}}(t)\Delta_{\R}(t)
\end{equation*}
Because the quotient representation $\Gamma \lra \Aut(\R)$ is trivial, this is simply the \emph{untwisted} Alexander polynomial associated to a once-punctured torus bundle.  Because $M$ is a 3-manifold with whose untwisted homology group $H^{1}(M,\Z)$ is rank 1, then $1$ is not a root of this untwisted polynomial and the untwisted polynomial itself is given by $k/(1-t)$ where $k$ is order of torsion of the manifold, namely $|\text{Tor } H^{1}(M;\Z)|$.    \\
\\
In light of the above remark, this means the twisted Alexander polynomial in $\mathfrak{gl}(4,\R)$ has one less factor of $(1-t)$ than the twisted polynomial in $\mathfrak{sl}(4,\R)$.  Moreover, this polynomial will be same as the one in $V\otimes V$ as $V$ and $V^{*}$ are isomorphic.  Calculating the twisted Alexander polynomial of $\Gamma \lra \Aut(V\otimes V)$ is straightforward in practice if we have the holonomy representation $\rho: \Gamma \lra \SO(3,1)$, as the matrix representation of $\rho(g)\otimes \rho(g) \in \Aut(V\otimes V)$ is simply given by the \emph{Kronecker Product} of $\rho(g)$ with itself.  By the above remarks and Theorem \ref{thm2} infinitesimal projective rigidity is guaranteed if the twisted Alexander polynomial in $\mathfrak{gl}(4,\R)$ is degree $5-1 = 4$.  We summarize these calculations in the following theorem.
\begin{theorem}\label{thm3}
Let $M$ be a hyperbolic once-punctured torus bundle with holonomy $\rho: \Gamma \lra \emph{\PSO}(3,1) \subset \emph{\PGL}(4,\R)$.  Lift the representation $\rho$ to $\tilde{\rho}: \Gamma \lra \SO(3,1)$.  If the twisted Alexander polynomial of $\tilde{\rho}\otimes \tilde{\rho}$ has $1$ as a root with multiplicity 4, then $M$ is infinitesimally projectively rigid rel cusp.  
\end{theorem}
\end{subsection}
\end{section}

\begin{section}{Calculations and Examples}\label{secCalcandExs}
In this section we employ Theorem \ref{thm2} and Theorem \ref{thm3} to a couple of hyperbolic once-punctured torus bundles to verify infinitesimal projective rigidity.  These calculations are done with the assistance of Mathematica.  Files containing relevant code can be found at \href{https://www.math.brown.edu/cdaly2/research.html}{Daly}.

\begin{subsection}{$L^{2}R^{2}$}\label{ssL2R2}
Let $M$ be the hyperbolic once-punctured torus bundle whose holonomy in $\PSO(3,1)$ is calculated Example \ref{exLLRR}.  We consider the representation $\Ad \rho: \Gamma \lra \Aut(\slf(4,\R))$.  Calculating explicit matrices for $\Ad\rho$ on the generators $a,b$ and $x$ is a matter of solving a linear system in 15-variables corresponding to how conjugation by $\rho(g)$ acts on $\slf(4,\R)$ for $g = a,b$ and $x$.  These calculations are carried out in \texttt{L2R2}\underscore\texttt{example.nb}.  \\
\\
To calculate $H^{1}(F,\del F; \g)$ we calculate the nullspace of the matrix associated to the linear transformation in Equation \ref{eqresl}, namely the kernel of the restriction map $\res_{l}^{B} : Z^{1}(F_{2};\g) \lra Z^{1}(\Z; \g)$.  This is an invariant subspace under the action induced by the monodromy of $\phi$.  Because we wish to use Theorem \ref{thm2} and calculate the multiplicity of 1 as a root of the characteristic polynomial of the $\phi$-action, we may instead look at the action of $\phi^{-1}$.  This is computationally a bit easier insofar as we don't have to take Fox-derivatives of as many inverses.  \\
\\
Explicitly, if $f$ is a cocycle in $H^{1}(F,\del F; \g)$, then $\phi^{-1}f$ acts via $(\phi^{-1}f)(a) = x^{-1}f(\phi(a))$.  Calculating this value is then a matter of taking Fox-derivatives.  These calculations are carried out below where $\phi(a) = xax^{-1} = ab^{2}$ and $\phi(b) = xbx^{-1} = bab^{2}ab^{2}$.  
\begin{align*}
\begin{array}{ll}
\frac{\del}{\del a}\phi(a) = 1& \frac{\del}{\del a}\phi(b) = a + ab \\
\phantom{==} & \phantom{==} \\
\frac{\del}{\del a}\phi(b) = b+bab^{2} & \frac{\del}{\del b}\phi(b) = 1 + ba+bab+bab^{2}a + bab^{2}ab
\end{array}
\end{align*}
So $\phi^{-1}$ takes an element $(X,Y) \in H^{1}(F,\del F; \g)$ and sends it to
\begin{equation*}
\phi^{-1}(X,Y) = \left(x^{-1}\left(X + \left(a+ab\right)Y\right), x^{-1}\left(\left(b+bab^{2}\right)X + \left(1+ba+bab+bab^{2}a+bab^{2}ab\right)Y\right)\right)
\end{equation*}
The characteristic polynomial of $\phi^{-1}$ is calculated in \texttt{L2R2.nb}.  It is given by
\begin{equation*}
p_{L^{2}R^{2}}(t) = -(t-1)^5 \left(t^2-18 t+1\right)^2 \left(t^6-114 t^5-17 t^4-316 t^3-17 t^2-114 t+1\right)
\end{equation*}
Thus by Theorem \ref{thm2} the once-punctured hyperbolic 3-manifold with monodromy $L^{2}R^{2}$ is infinitesimally projectively rigid rel cusp.
\end{subsection}

\begin{subsection}{$R^{2}L$}\label{ssR2L}
Let $M$ be the hyperbolic once-punctured torus bundle with monodromy $R^{2}L$ where $R$ and $L$ are defined in Example \ref{exLLRR}.  The fundamental group $\Gamma$ admits a presentation of the form $\Gamma = \langle a,b,x \, | \, xax^{-1} = aba^{2} \text{ and } xbx^{-1} = ba^{2}\rangle$.  Using the trace-relations we see that
\begin{equation*}
(AC - B)A - C = A \phantom{==} (AC-B) = B \phantom{==}A^{2} + B^{2} + C^{2} = ABC
\end{equation*}
where $A,B,C$ are the traces of $a,b$ and $ab$ respectively.  Solving this system in $A,B$ and $C$ and using the Maclachlan-Reid free group representation yields
\begin{small}
\begin{equation*}
\rho(a) = \left(
\begin{array}{cc}
 -\frac{1}{2} \sqrt{\frac{1}{2} \left(5-i \sqrt{7}\right)} & \frac{6+2 i \sqrt{7}}{8
   \sqrt{10-2 i \sqrt{7}}} \\
 -\sqrt{\frac{1}{2} \left(5-i \sqrt{7}\right)} & -\frac{1}{2} \sqrt{\frac{1}{2} \left(5-i
   \sqrt{7}\right)} \\
\end{array}
\right),
\rho(b) = \left(
\begin{array}{cc}
 \frac{\sqrt{7}+13 i}{2 \sqrt{7}+2 i} & \frac{1}{8} i \left(\sqrt{7}+i\right) \\
 -\frac{\sqrt{7}+5 i}{\sqrt{7}+i} & \frac{1}{4} \left(1-i \sqrt{7}\right) \\
\end{array}
\right)
\end{equation*}
\end{small}
Solving for the meridian imposed by the equations $xax^{-1} = aba^{2}$ and $xbx^{-1} = ba^{2}$ yields that
\begin{equation*}\rho(x) = 
\left(
\begin{array}{cc}
 -1 & -\frac{1}{4} i \left(\sqrt{7}-i\right) \\
 0 & -1 \\
\end{array}
\right)
\end{equation*}
Applying the isomorphism between $\PSL(2,\C)$ and $\PSO(3,1)$ yields the holonomy representation into $\PSO(3,1)$ below.  
\begin{scriptsize}
\begin{equation}\label{eqsor2l}
\rho(a) = \left(
\begin{array}{cccc}
 \frac{5}{4 \sqrt{2}} & \frac{\sqrt{\frac{7}{2}}}{4} & \frac{\sqrt{\frac{7}{2}}}{8} &
   \frac{\sqrt{\frac{7}{2}}}{8} \\
 \frac{\sqrt{\frac{7}{2}}}{4} & \frac{3}{4 \sqrt{2}} & -\frac{17}{8 \sqrt{2}} &
   \frac{15}{8 \sqrt{2}} \\
 -\frac{\sqrt{\frac{7}{2}}}{8} & \frac{17}{8 \sqrt{2}} & -\frac{17}{16 \sqrt{2}} &
   \frac{31}{16 \sqrt{2}} \\
 \frac{\sqrt{\frac{7}{2}}}{8} & \frac{15}{8 \sqrt{2}} & -\frac{31}{16 \sqrt{2}} &
   \frac{49}{16 \sqrt{2}} \\
\end{array}
\right)
 \rho(b) = \left(
\begin{array}{cccc}
 -\frac{3}{4} & -\frac{\sqrt{7}}{4} & -\frac{\sqrt{7}}{2} & \frac{\sqrt{7}}{2} \\
 \frac{3 \sqrt{7}}{4} & -\frac{5}{4} & \frac{17}{4} & -\frac{19}{4} \\
 -\frac{\sqrt{7}}{4} & 0 & \frac{15}{16} & -\frac{9}{16} \\
 -\frac{3 \sqrt{7}}{4} & 1 & -\frac{71}{16} & \frac{81}{16} \\
\end{array}
\right)
\rho(x) = \left(
\begin{array}{cccc}
 1 & 0 & -\frac{\sqrt{7}}{4} & -\frac{\sqrt{7}}{4} \\
 0 & 1 & \frac{1}{4} & \frac{1}{4} \\
 \frac{\sqrt{7}}{4} & -\frac{1}{4} & \frac{3}{4} & -\frac{1}{4} \\
 -\frac{\sqrt{7}}{4} & \frac{1}{4} & \frac{1}{4} & \frac{5}{4} \\
\end{array}
\right)
\end{equation}
\end{scriptsize}
We may calculate the Kronecker product of these matrices to yield the representation of $\Gamma$ into $V\otimes V$ where $V$ is the representation $\rho : \Gamma \lra \SO(3,1)$ given in Equation \ref{eqsor2l} above.  We calculate the twisted Alexander polynomial of this representation via Alexander matrix.  To simplify calculations we may instead use the Fox-derivatives of the relations $R_{1} = \phi(a)xa^{-1}x^{-1}$ and $R_{2} = \phi(b)xb^{-1}x^{-1}$.  This yields 
\begin{equation}\label{eqtwistedminorcalc}
\Delta_{\rho\otimes\rho,\alpha}(t) = \det \left( \begin{array}{ccc}
\frac{\del}{\del a} \phi(a) - tx & \frac{\del}{\del b} \phi(a)  \\
\frac{\del}{\del a}\phi(b)   &  \frac{\del}{\del b} \phi(b) - tx 
\end{array}\right)/\det(\text{Id}_{\mathfrak{gl}(4,\R)}-tx)
\end{equation}
where above we are suppressing the representation $\rho\otimes\rho$ to save space.  This calculation yields
\begin{small}
\begin{equation*}
\Delta_{\rho\otimes\rho,\alpha}(t)  =(t-1)^4 \left(t^2-4 t+1\right) \left(t^4-18 t^3+90 t^2-18 t+1\right) \left(t^6-38 t^5+15
   t^4-84 t^3+15 t^2-38 t+1\right)
\end{equation*}
\end{small}
Thus by Theorem \ref{thm3}, the hyperbolic once-punctured torus bundle with monodromy $R^{2}L$ is infinitesimally projectively rigid rel cusp.  If we carry out the calculations of the characteristic polynomial of $\phi$ acting on $H^{1}(F,\del F)$ similar to the case as we did in Section \ref{ssL2R2}, we obtain the twisted Alexander polynomial below.  
\begin{small}
\begin{equation*}
-(t-1)^5 \left(t^4-18 t^3+90 t^2-18 t+1\right) \left(t^6-38 t^5+15 t^4-84 t^3+15 t^2-38
   t+1\right)
\end{equation*}
\end{small}
Note it differs from the one on $\mathfrak{gl}(4,\R)$ by a factor of $t^2-4 t+1$, which, when evaluated at $t = 1$ equals 2.  This is the order of the homology group of $M$ with monodromy $R^{2}L = \left(
\begin{array}{cc}
 2 & 3 \\
 1 & 2 \\
\end{array}
\right)$.  This is because the holonomy has order $|\tr \phi - 2| = 2$ \cite{Miller2023Azumaya}.  These calculations are carried out in \texttt{R2L}\underscore\texttt{example.nb}
\end{subsection}

\end{section}

\begin{section}{Bibliography}
\bibliographystyle{plain}
\bibliography{refs}

\begin{thebibliography}{10}

\bibitem{1914Vorlesungen}
Vorlesungen über die theorie der automorphen funktionen.
\newblock {\em Monatshefte für Mathematik und Physik}, 25(1):A34--A34, 12
  1914.

\bibitem{Ballas2018Convex}
Samuel Ballas, Jeffrey Danciger, and Gye-Seon Lee.
\newblock Convex projective structures on nonhyperbolic three-manifolds.
\newblock {\em Geometry \&; Topology}, 22(3):1593--1646, 3 2018.

\bibitem{Ballas2014Deformations}
Samuel~A Ballas.
\newblock Deformations of noncompact projective manifolds.
\newblock {\em Algebraic \&; Geometric Topology}, 14(5):2595--2625, 11 2014.

\bibitem{BODEN_FRIEDL_2014}
Hans~U. Boden and Stefan Friedl.
\newblock Metabelian sl(n, $\mathbb c$) representations of knot groups iv:
  twisted alexander polynomials.
\newblock {\em Mathematical Proceedings of the Cambridge Philosophical
  Society}, 156(1):81–97, 2014.

\bibitem{Brown1982Cohomology}
Kenneth~S. Brown.
\newblock {\em Cohomology of Groups}.
\newblock Springer New York, 1982.

\bibitem{Cooper2006Computing}
Daryl Cooper, Darren Long, and Morwen Thistlethwaite.
\newblock Computing varieties of representations of hyperbolic 3-manifolds into
  sl(4, r).
\newblock {\em Experimental Mathematics}, 15(3):291--305, 1 2006.

\bibitem{CULLER198664}
Marc Culler.
\newblock Lifting representations to covering groups.
\newblock {\em Advances in Mathematics}, 59(1):64--70, 1986.

\bibitem{Daly}
Charles Daly.
\newblock Prfe files, 2024.
\newblock \url{https://arxiv.org/abs/2408.08405} [Accessed: (Oct 2024)].

\bibitem{BDLM}
Jeff Danciger, Sam Ballas, Gye-Seon Lee, and Ludovic Marquis.
\newblock Properly convex dehn filling, 2024.
\newblock \url{https://www.math.fsu.edu/~ballas/} [Accessed: October 2024].

\bibitem{Dunfield2012Twisted}
Nathan~M. Dunfield, Stefan Friedl, and Nicholas Jackson.
\newblock Twisted alexander polynomials of hyperbolic knots.
\newblock {\em Experimental Mathematics}, 21(4):329--352, 12 2012.

\bibitem{Fox1953Free}
Ralph~H. Fox.
\newblock Free differential calculus. i: Derivation in the free group ring.
\newblock {\em The Annals of Mathematics}, 57(3):547, 5 1953.

\bibitem{Goldman2020Parallelism}
William Goldman.
\newblock Parallelism on lie groups and fox’s free differential calculus,
  2020.

\bibitem{Goldman2009Trace}
William~M. Goldman.
\newblock Trace coordinates on fricke spaces of some simple hyperbolic
  surfaces.
\newblock In {\em IRMA Lectures in Mathematics and Theoretical Physics}, pages
  611--684. EMS Press, 3 2009.

\bibitem{Hatcher2002Algebraic}
Allen Hatcher.
\newblock {\em Algebraic Topology}.
\newblock Cambridge University Press, 2002.

\bibitem{gomtrie2016SomeRR}
Michael Heusener.
\newblock Some recent results about the sln(c)–representation spaces of knot
  groups.
\newblock 2016.

\bibitem{Heusener2011Infinitesimal}
Michael Heusener and Joan Porti.
\newblock Infinitesimal projective rigidity under dehn filling.
\newblock {\em Geometry \&; Topology}, 15(4):2017--2071, 10 2011.

\bibitem{Hodgson}
Craig Hodgson.
\newblock {\em Degeneration and Regeneration of Geometric Structures on
  Three-Manifolds}.
\newblock Phd thesis, Princeton, Princeton, New Jersey, June 1986.
\newblock Available at
  \url{http://homepages.warwick.ac.uk/~masgar/Maths/hodgson.pdf}.

\bibitem{Mathematica}
Wolfram~Research{,} Inc.
\newblock Mathematica, {V}ersion 14.0.
\newblock Champaign, IL, 2024.

\bibitem{Johnson1987Deformation}
Dennis Johnson and John~J. Millson.
\newblock Deformation spaces associated to compact hyperbolic manifolds.
\newblock In {\em Discrete Groups in Geometry and Analysis}, pages 48--106.
  Birkhäuser Boston, 1987.

\bibitem{Jorgensen}
Troels J{\o}rgensen.
\newblock {\em On pairs of once-punctured tori}, page 183–208.
\newblock London Mathematical Society Lecture Note Series. Cambridge University
  Press, 2003.

\bibitem{Kapovich2010Hyperbolic}
Michael Kapovich.
\newblock {\em Hyperbolic Manifolds and Discrete Groups}.
\newblock Birkhäuser Boston, 2010.

\bibitem{Kitano}
Teruaki Kitano.
\newblock Winter braids lecture notes.
\newblock \url{https://arxiv.org/pdf/1510.03216}, 2016.

\bibitem{Lin2001RepresentationsOK}
Xiaoxia Lin.
\newblock Representations of knot groups and twisted alexander polynomials.
\newblock {\em Acta Mathematica Sinica}, 17:361--380, 2001.

\bibitem{claraloh}
Clara L\'{o}h.
\newblock Group cohomology.
\newblock
  \url{https://loeh.app.uni-regensburg.de/teaching/grouphom_ss19/lecture_notes.pdf},
  2019.

\bibitem{Maclachlan2003Arithmetic}
Colin Maclachlan and Alan~W. Reid.
\newblock {\em The Arithmetic of Hyperbolic 3-Manifolds}.
\newblock Springer New York, 2003.

\bibitem{Miller2023Azumaya}
Nicholas Miller.
\newblock Azumaya algebras and once-punctured torus bundles.
\newblock 2023.

\bibitem{sikora}
Adam~S. Sikora.
\newblock Character varieties.
\newblock {\em arXiv}, 2009.

\bibitem{Thurston2022Geometry}
William~P. Thurston.
\newblock {\em The Geometry and Topology of Three-Manifolds}.
\newblock American Mathematical Society, 7 2022.

\bibitem{WADA1994241}
Masaaki Wada.
\newblock Twisted alexander polynomial for finitely presentable groups.
\newblock {\em Topology}, 33(2):241--256, 1994.

\bibitem{Weil1964Remarks}
Andre Weil.
\newblock Remarks on the cohomology of groups.
\newblock {\em The Annals of Mathematics}, 80(1):149, 7 1964.

\end{thebibliography}
\end{section}
\end{document}